\documentclass[12pt]{article}

\usepackage[letterpaper,margin =1in]{geometry}

\usepackage{amssymb}
\usepackage{amsmath}
\usepackage{amsrefs}
\usepackage{amsthm}

\usepackage{mathtools}
\usepackage{enumerate}
\usepackage{fancyhdr}
\usepackage{braket}

\usepackage{hyperref}
\usepackage{titlesec}
\usepackage{xcolor}

\DeclareFontFamily{OT1}{pzc}{}
\DeclareFontShape{OT1}{pzc}{m}{it}%
             {<-> s * [1.200] pzcmi7t}{}
      
\newlength\tindent
 \newcommand{\bigO}{\mathcal{O}}
\newcommand\eIL[1]{{\color{black}#1}}

\DeclareMathAlphabet{\mathscr}{OT1}{pzc}{m}{it}

\providecommand{\eps}{\varepsilon}

\newenvironment{frcseries}{\fontfamily{frc}\selectfont}{}

\theoremstyle{remark}
\newtheorem{definition}{Definition}
\newtheorem{remark}{Remark}
\newtheorem{assumption}{Assumption}
\theoremstyle{plain}
\newtheorem{lemma}{Lemma}
\newtheorem{theorem}{Theorem}

\begin{document}

\thispagestyle{plain}

{\Large Contact singularities in nonstandard slow-fast dynamical systems}\\

\begin{center}
Ian Lizarraga, Robby Marangell, and Martin Wechselberger\\
{\it\small School of Mathematics and Statistics, University of Sydney, Camperdown 2006, Australia}
\end{center}

\begin{abstract}
We develop the contact singularity theory for singularly-perturbed (or `slow-fast') vector fields of the general form $z' = H(z,\eps)$, $z\in\mathbb{R}^n$ and $\eps\ll 1$. Our main result is the derivation of computable, coordinate-independent defining equations for contact singularities under an assumption that the leading-order term of the vector field admits a suitable factorization. This factorization can in turn be computed explicitly in a wide variety of applications. We demonstrate these computable criteria by locating contact folds and, for the first time, contact cusps in some nonstandard models of biochemical oscillators. 	
\end{abstract}

\section{Introduction} 


Classifying the loss of normal hyperbolicity of the critical manifold is a fundamental step in the analysis of 
slow-fast dynamical systems. For systems in the so-called {\it standard} form\footnote{Throughout the paper we use prime notation $'=d/dt$ to denote derivates with respect to the (fast) time variable $t$ and the notation $D_{v}$ to denote partial derivatives with respect to phase space variables $v$.}
\begin{eqnarray}
x' &=&  \eps g(x,y,\eps) \nonumber \\
y' &=& f(x,y,\eps), \label{eq:standard}
\end{eqnarray}

 the $k$-dimensional {\it critical manifold} lies inside the zero set of a smooth mapping $f(x,y,0): \mathbb{R}^{n}\to \mathbb{R}^{n-k}$.  Loss of normal hyperbolicity occurs along points where the critical manifold becomes tangent to the {\it layer problem} of \eqref{eq:standard}, formally defined by the $\eps \to 0$ limit:
\begin{eqnarray}
x' &=& 0  \nonumber\\
y' &=& f(x,y,0). \label{eq:layerstandard}
\end{eqnarray}

In the planar case $x,y\in \mathbb{R}$, the solutions of the layer flow consist of trajectories lying within vertical lines $x = c$ (also known as {\em fast fibers}). In general, solutions of the layer flow lie in hyperplanes orthogonal to the coordinate axes of the slow variable $x\in \mathbb{R}^k$. \\

 Of particular interest  is the loss of normal hyperbolicity associated with a geometric `fold' structure in phase space where the layer flow has a  tangency with the critical manifold, which allows switching between slow and fast motion as observed in, e.g., relaxation oscillations.  The defining equations for an isolated fold point $F\in S$ in the planar, standard slow-fast system \eqref{eq:standard} are well-known (see for eg. Sec. 8.1 in \cite{kuehn2015}):
 \begin{equation} \label{eq:foldconds}
\begin{aligned}
f|_F= 0\,,  & \qquad D_yf|_F = 0,\\ 
D_{yy}f|_F\, \neq 0, & \qquad D_xf|_F \neq  0.\nonumber 
\end{aligned} 
\end{equation}

The geometric content of the first three conditions (on the derivatives with respect to $y$) is that the critical manifold of equilibria lying inside the zero set $\{f(x,y,0) = 0\}$ makes parabolic contact with the vertical layer flow at $F$, whereas the role of the final transversality condition can be deduced from the classical singularity theory of the generic fold map $f(\alpha,x) = \alpha + x^2$ (see for eg. \cite{kuznetsov}): the slow variable $x$ plays the role of the unfolding parameter.\\
 
Takens  \cite{takens1975,takens1976a,takens1976b} began the classification of singularities of low codimension from the point of view of constrained differential equations, corresponding to the $\eps = 0$ limit of the equivalent slow standard system
 \begin{equation}  \label{eq:standard-slow}
 \begin{aligned}
\dot{x} &=&  g(x,y,\eps) \\
\eps \dot{y} &=& f(x,y,\eps),
\end{aligned}
\end{equation}
 known as the {\em reduced problem} in the geometric singular perturbation theory (GSPT) literature. Here, the dot notation $\,\dot{}\, = d/d\tau$ denotes the derivative with respect to the slow time $\tau = \eps t$. Since then, much work has been done to extend this local analysis for $0 < \eps \ll 1$,  in the case where the slow variables $x \in \mathbb{R}^{k}$ play the role of unfolding parameters for folds \cite{benoit1983,krupa2001,szmolyan2001,szmolyan2004,wechselberger2012} and cusps \cite{broer2013,jardon2016}. In a complementary direction, the theory of {\it bifurcations without parameters} \cite{liebscher2015} has recently been developed to describe classical bifurcations in terms of breakdown of normal hyperbolicity of manifolds of equilibria. \\

The purpose of this work is to provide a more general classification of loss of normal hyperbolicity of the critical manifold for the larger class
\begin{eqnarray}
z' &=& H(z,\eps) \label{eq:nonstandard}
\end{eqnarray} 
of {\it nonstandard} multiple-timescale systems, by making use of a suitably general analogue of the classical singularity theory. The system \eqref{eq:nonstandard} defines a singular perturbation problem if the set of equilibria of the corresponding {\it layer problem} 
\begin{eqnarray}
z' &=& H(z,0) \label{eq:layer}
\end{eqnarray}
contains a differentiable manifold $S$, which we refer to as the critical manifold of \eqref{eq:nonstandard}. The relationship between Eqs. \eqref{eq:standard} and \eqref{eq:nonstandard} is that coordinate transformations placing \eqref{eq:nonstandard} in the form \eqref{eq:standard} are defined only locally; in other words, there is in general no globally defined coordinate splitting into `slow' versus `fast' directions. As in the standard case, solutions of the layer problem determine the leading order fast motion relative to $S$; however, these solutions are no longer `unidirectional'; i.e., they are no longer restricted to lie in hyperplanes orthogonal to a set of distinguished `slow' coordinate axes, but rather are allowed to bend and curve throughout the phase space. \\

 To illuminate the new complication, consider the following planar system:
\begin{eqnarray}
\begin{pmatrix} x' \\ y' \end{pmatrix} &=& \begin{pmatrix}  x^3-2x^2+xy-x-2y+2 \\ y+x^2-1 \end{pmatrix} +  \eps \begin{pmatrix} G_1(x,y,\eps) \\ G_2(x,y,\eps) \end{pmatrix}  \label{eq:parab}
\end{eqnarray}
where $G_1,G_2$ are smooth functions in their arguments. System \eqref{eq:parab} is of the form \eqref{eq:nonstandard} with $z = (x,y)$.
The corresponding layer problem \eqref{eq:layer} can be factorised:
\begin{eqnarray}
\begin{pmatrix} x' \\ y' \end{pmatrix} &=&  \begin{pmatrix}  x^3-2x^2+xy-x-2y+2 \\ y+x^2-1 \end{pmatrix} \label{eq:parablayer}\\
&=&\begin{pmatrix} (y+x^2-1)(x-2) \\ y+x^2-1 \end{pmatrix} \nonumber\\
&=& \begin{pmatrix} x-2 \\ 1 \end{pmatrix} (y+x^2-1) = N(x,y)f(x,y). \nonumber
\end{eqnarray}

Important geometric content about the solutions of \eqref{eq:parablayer} is encoded in this factorisation. For example, since $N(x,y)$ is everywhere nonzero, the critical manifold can be read off from $f(x,y)$:
\begin{eqnarray}
S = \{(x,y) \in \mathbb{R}^2:  f(x,y) = 0\} = \{(x,y)\in \mathbb{R}^2: y  = 1-x^2\}. \label{eq:parabcrit}
\end{eqnarray}

Away from $S$, we can formally rescale time by the scalar function $f(x,y)$ in \eqref{eq:parablayer} to obtain a desingularised layer problem:
 \begin{eqnarray}
\begin{pmatrix} x' \\ y' \end{pmatrix} &=& N(x,y) = \begin{pmatrix}  x-2 \\ 1 \end{pmatrix}. \label{eq:parabdesing}
\end{eqnarray}

\begin{remark}
The supports of \eqref{eq:parablayer} and \eqref{eq:parabdesing} are identical, but the solution curves of \eqref{eq:parabdesing} are regular owing to the removal of the singularity set $\{f = 0\}$. Away from this set, the trajectories of the two systems are identical up to time rescalings and possible orientation reversals due to sign changes of $f$.
\end{remark}

\begin{figure}[t]
  \centering
        \includegraphics[height=0.55\textwidth,width=0.95\textwidth]{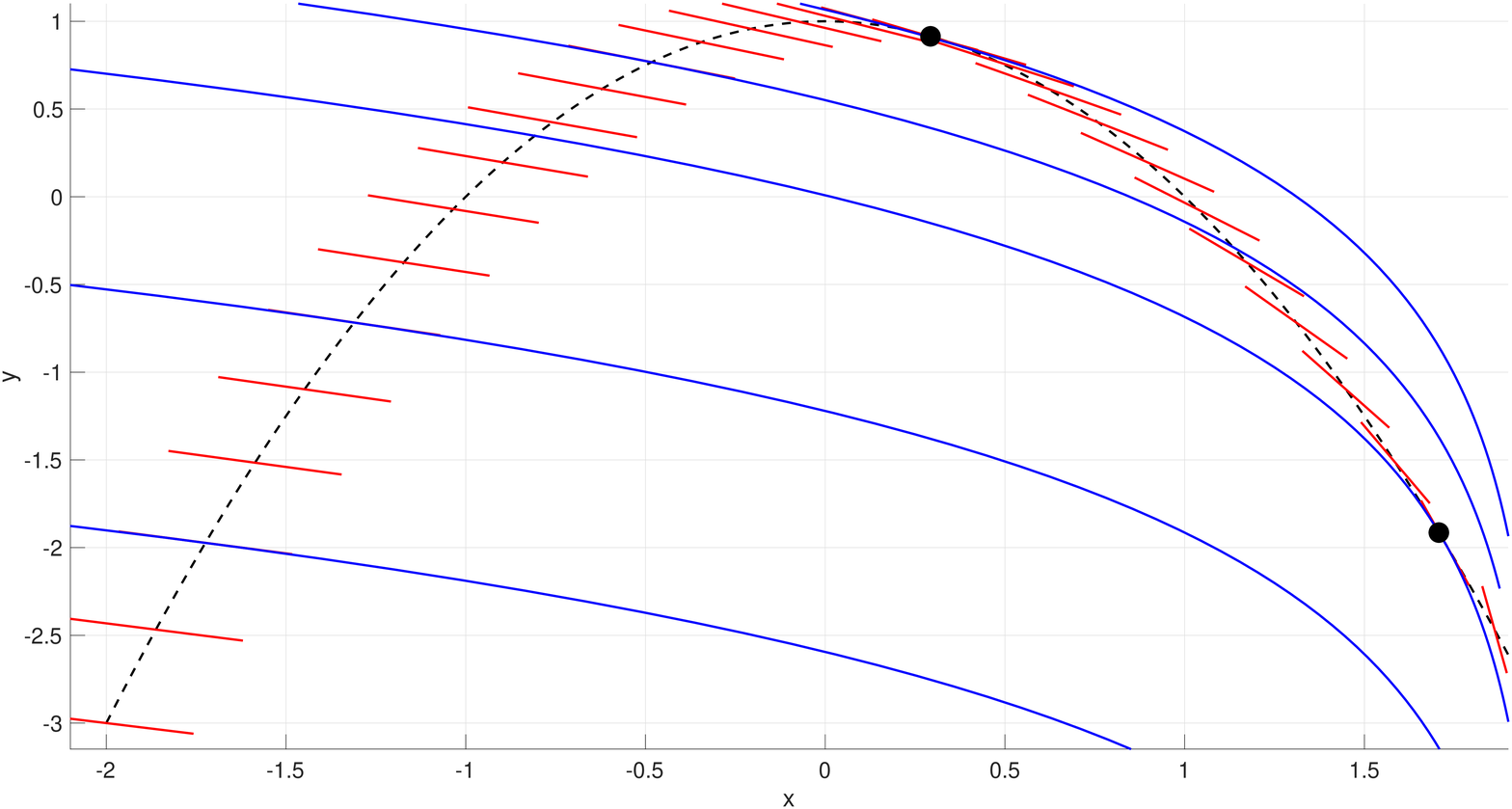}
      \caption{Critical manifold (dashed black curve) $S = \{(x,y): y = 1-x^2\}$ of system \eqref{eq:parab} together with trajectories of the desingularised layer problem \eqref{eq:parabdesing} (blue curves). Red line segments denote portions of the tangent spaces of the desingularised layer trajectories at basepoints along $S$. Tangencies of the layer solutions with the critical manifold are labeled with black points.}  
      \label{fig:parab}
\end{figure}

At each point $(x,y) \in S$, the tangent space of a (regular) solution curve of \eqref{eq:parabdesing} at $(x,y)$ is spanned by the nonzero vector $N(x,y)$. The corresponding time-rescaled solutions of \eqref{eq:parablayer} must therefore approach $S$ tangent to $N(x,y)$. In Fig. \ref{fig:parab} we overlay the critical manifold with solution curves of \eqref{eq:parabdesing}. 

\begin{remark}
The fast fibers of \eqref{eq:parabdesing} can be calculated explicitly here and are given by a family of curves $y = \ln|x-2| + C$ (for $x \neq 2$) and $x = 2$.
\end{remark}

The novelty is that the geometric `fold' point of the critical manifold \eqref{eq:parabcrit} at $(0,1)$ does {\it not} correspond to a tangency of the critical manifold with the layer flow! 
In contrast to the standard theory, there is no longer any correlation between geometric folds of the critical manifold and tangencies with the layer flow.

On the other hand, there appears to be `fold-like' behavior with $S$ at two points. 
Evidently, the layer solutions approach $S$ in a tangent direction precisely when $N$ lies in the kernel of $D\!f  $:
\begin{eqnarray*}
\left. D\!f   N \right\vert_S &=& 0\\
2x(x-2) + 1 &=& 0\\
\Rightarrow x &=& \frac{1}{2}(2 \pm \sqrt{2}).
\end{eqnarray*}

The geometric constraint $D\!f N|_S = 0$ immediately generalises one of the classical fold conditions 
\begin{eqnarray*}
D_yf(0,0,0) &=& 0
\end{eqnarray*}
in \eqref{eq:foldconds}. Indeed, note that $N = \begin{pmatrix} 0 \\ 1 \end{pmatrix}$ for planar systems in the standard form, and thus
\begin{eqnarray*}
D\!f   N &=& \begin{pmatrix} D_x f & D_y f \end{pmatrix} \begin{pmatrix} 0 \\ 1 \end{pmatrix} = D_y f.
\end{eqnarray*}

This computation suggests that the analogous second-order nondegeneracy criterion generalizing 
\begin{eqnarray*}
D_{yy}f(0,0,0) &\neq&0
\end{eqnarray*}
 will likely involve second derivatives of $f$, measuring curvature of the critical manifold, {\it as well as} derivatives of $N$, measuring the curvature of the fibers. These extra contributions make higher-order degeneracies harder to identify. 
 
The fold points identified in Fig. \ref{fig:parab} are {\it contact points of order one}, for which defining equations and genericity conditions have recently been computed  \cite{wechselberger2020}.  These results explicitly use the factorisation 
\begin{eqnarray}
z' &=& N(z)f(z) \label{eq:factorintro}
\end{eqnarray}

of the layer problem. Such factorisations are natural in many applications, including chemical reaction networks \cite{schauer1983, stiefenhofer}. Motivated by these reaction network models, Goeke and Walcher \cite{goeke2014} have recently adapted algorithms from computational algebraic geometry to explicitly construct such factorisations (locally) in the case of rational vector fields.  \\

 The present paper has two goals. The first is to develop the contact singularity theory for slow-fast systems, giving a firm theoretical foundation for the computations in \cite{wechselberger2020}.  Mather introduced the general notion of contact singularities between equidimensional manifolds in \cite{mather1969b}, but here we adopt the extended development of Montaldi's PhD thesis \cite{montaldi1983,montaldi1986,montaldi1991}.  The basic objects of study in Montaldi's framework are {\it contact maps} between smooth manifolds of not necessarily equal dimension. A key result of the present paper is the establishment of a rigorous relationship between these computations and the singularity classes of the contact map between the fast fibers and the slow manifold. In fact, for the case considered in this paper where the fast fibers make contact with the critical manifold along a (one-dimensional) curve, the singularity classes turn out to be the well-known $A_c$ singularities \cite{guseinzade}. The nondegeneracy condition for the $A_c$ singularity of the contact map is then identical to the corresponding codimension-$k$ folded singularity nondegeneracy condition, up to a trivial projection along the center flow.
 
 \begin{remark}
These are typically referred to as ``$A_k$'' singularities in the literature, but we use the variable $c$ to avoid confusion with the dimension $k$ of the critical manifold.\\
\end{remark}

 \eIL{The second goal is to extend the results for defining equations of the contact fold in \cite{wechselberger2020}. In the present paper we give coordinate-independent defining equations for slow unfoldings of contact points of arbitrary order.  We also identify {\it contact cusps} in two nonstandard models of biochemical oscillators.} \\

The present paper proceeds as follows: In Section \ref{sec:mtds}, we give an account of multiple-timescale systems having a nonstandard slow-fast splitting. In Section \ref{sec:contact}, we give a rigorous definition of contact between a one-dimensional manifold and a $k$-dimensional manifold in $\mathbb{R}^n$, where $1 \leq k < n$, culminating in a full description of both the singularity classes and computable defining equations for contact folds and cusps. In Section \ref{sec:examples}, we locate contact folds and cusps in several three-dimensional examples. We first test our results on standard slow-fast systems, and then we consider two nonstandard models of biochemical oscillators with negative feedback loops: a minimal three-component model, and a mitotic oscillator model. In both cases, we demonstrate the existence of contact cusps previously not shown.
We conclude in Sec. \ref{sec:conclusion}.

\section{Multiple-timescale dynamical systems} \label{sec:mtds}
\subsection{The nonstandard formulation}
We begin by giving an abbreviated treatment of nonstandard multiple-timescale dynamical systems, following Fenichel's seminal work on geometric singular perturbation theory \cite{fenichel1979}, and Wechselberger's more recent treatment \cite{wechselberger2020} which extends the framework to loss of normal hyperbolicity. Consider the family of vector fields \eqref{eq:nonstandard}, formally expanded in $\eps$:
\begin{eqnarray}
z' &=& h(z) + \eps G(z,\eps). \label{eq1}
\end{eqnarray}

\begin{definition} The family of vector fields \eqref{eq1} is called a {\it singular perturbation problem} if the set of equilibria $\{z \in \mathbb{R}^n: h(z) = 0\}$ contains a $k$-dimensional differentiable manifold for some $1 \leq k < n$. 
\end{definition}

Our first assumption restricts the geometry of this equilibrium set.

\begin{assumption} \label{ass1} The system \eqref{eq1} is a singular perturbation problem with a single subset $S \subset \{z \in \mathbb{R}^n: h(z) = 0\}$ that forms a connected $k$-dimensional differentiable manifold, called the {\it critical manifold}. \end{assumption}

\begin{definition} \label{def:layer}
Given a singular perturbation problem \eqref{eq1}, the corresponding singular limit system
\begin{eqnarray} 
z' &=& H(z,0) = h(z) \label{layer}
\end{eqnarray}
is called the {\it layer problem} of \eqref{eq1}.
\end{definition}

For convenience we distinguish two subsets in the spectrum of $Dh$.
\begin{definition} \label{def:eigs}
The {\it trivial eigenvalues} of the layer problem at $z\in S$ are the $k$ zero eigenvalues of $Dh(z)$  corresponding to the dimension of the tangent space $T_zS$. The remaining $n-k$ eigenvalues along $S$ are called {\it nontrivial}.
\end{definition}

Note that a nontrivial eigenvalue can be equal to zero.

\begin{definition} \label{def:contact} The set $S_n \subset S$ denotes the subset where all nontrivial eigenvalues of $Dh$ evaluated along $S_n$ are nonzero. \end{definition}
%
%

Along the set $S_n$, we may construct a pointwise-defined splitting of the tangent bundle along $z \in S_n$:
\begin{eqnarray}
T_z\mathbb{R}^n &=& T_z S \oplus N_z, \label{eq:splitting}
\end{eqnarray}

where $N_z$ is called the {\it linear fast fiber} at basepoint $z \in S_n$ identified with the quotient space $T_z\mathbb{R}^n/T_z S$. We may define the {\it tangent bundle of $S$} via the construction
\begin{eqnarray*}
TS &=& \cup_{z\in S_n} T_z S.
\end{eqnarray*}
The corresponding bundle 
\begin{eqnarray*}
N &=& \cup_{z\in S_n} N_z
\end{eqnarray*}

is called the {\it (linear) fast fiber bundle}. \\

Along $S_n$, we may define a projection operator 
\begin{eqnarray*}
\Pi^S: TS \oplus N \to TS.
\end{eqnarray*}

Given a point $z \in S_n$, the map $\Pi^S|_{z\in S_n}$ can be characterised geometrically as an oblique projection onto $T_z S$ along parallel translates of the fast fiber $N_z$ (see \cite{wechselberger2020,lizarraga2020}). \\

The layer flow near $S$ defines a locally invariant fast foliation $\mathscr{F}$ of the layer problem in a tubular neighborhood of $S$. The neighbourhood near a point of a $k$-dimensional critical manifold can be characterised geometrically as being foliated by regular level sets of $k$ smooth functions (see \cite{wechselberger2020} for details). In the case $k = n-1$, the tubular neighbourhood is foliated by the regular solution curves of a desingularised layer flow. We discuss this further in Sec. \ref{sec:straighten}.

Next, we assume a factorization of $h$ which captures the essential geometries of the slow and fast structures near to the critical manifold:

\begin{assumption} \label{ass3} The function $h(z)$ can be factorised as follows:
\begin{eqnarray}
h(z) &=& N(z)f(z), \label{factor}
\end{eqnarray}

where the $i$th column of the $n \times (n-k)$ matrix function $N(z)$,  $N_i = (N_i^1 ~ \cdots~N_{i}^{n})^T$ consists of smooth functions $N_i^k: \mathbb{R}^n \to \mathbb{R}$. Assume that $N(z)$ has full column rank $n-k$ for each $z \in S$, and furthermore that singularities of $N(z)f(z)$ for $z \notin S$ are isolated, if they exist. We furthermore assume that the critical manifold $S$ is  equal to the zero level set of a submersion $f: \mathbb{R}^n \to \mathbb{R}^{n-k}$:
\begin{eqnarray*}
S &=& f^{-1}(0).
\end{eqnarray*}
 \end{assumption}

\begin{remark} The question of existence and uniqueness of factorizations of the form \eqref{factor} for a singularly perturbed system of the general form \eqref{eq1} has only partial answers. Local factorizations can be constructed explicitly in the case the $h(z)$ is a rational vector field in $z$ \cite{goeke2014}. This includes a large variety of applied problems-- notably, many chemical reaction networks can be modeled in this framework \cite{schauer1983, stiefenhofer}. In practice, these local factorizations can often be shown {\it a posteriori} to hold over large open sets of the phase space. \end{remark}

The following two results immediately demonstrate the usefulness of this factorization:

\begin{lemma} \label{lemma:Nf1} For $z \in S_n$, the column vectors of $N(z)$ form a basis for the range of $Dh(z)$ and the transposes of the row vectors of $D\!f  (z)$ form a basis of the orthogonal complement of the kernel of $Dh(z)$ (i.e. a basis of the orthogonal complement of the tangent space $T_z S$). \end{lemma}

\begin{lemma}  \label{lemma:Nf2} The nontrivial eigenvalues of the layer problem of \eqref{eq1} along $S$ are equal as a set to the eigenvalues of $D\!f   N|_S$. \end{lemma}

Proofs. See \cite{wechselberger2020,lizarraga2020}. 

\subsection{The contact set}

 We are concerned with studying the subset  $S-S_n$, where the local tangent splittings \eqref{eq:splitting} break down due to the alignment of the fast fiber bundle with the critical manifold along a one-dimensional subspace. 
\begin{definition} \label{def:contactset}
The {\it contact set} $F \subset S-S_n$ is the set of points $z_0 \in S-S_n$ where exactly one nontrivial real eigenvalue of $Dh(z_0)$ (see Def. \ref{def:eigs}) vanishes. \end{definition}

From Lemmas \ref{lemma:Nf1}--\ref{lemma:Nf2} we have the following characterization of the contact set: 
\begin{eqnarray}
F &=& \{z \in S: \text{rank}(D\!f   N) = n - k - 1\}. \label{contactset}
\end{eqnarray}

Since $D\!fN$ is a matrix of size $(n-k)\times(n-k)$, a necessary condition is
\begin{eqnarray}
\det (D\!fN)|_{S} &=& 0. \label{eq:tangcond}
\end{eqnarray}

\begin{assumption} \label{ass2} The contact set $F$ is nonempty. \end{assumption}


Geometrically, $F$ is the set of points where a one-dimensional subspace of the fast fiber bundle locally aligns with the tangent space of $S$. In this setting it is straightforward to deduce the direction of tangency. 

\begin{lemma}
For $z_0 \in F$, let $r \in \mathbb{R}^{n-k}$ be any nontrivial column of adj$(D\!f(z_0) N(z_0))$, the adjugate of $D\!f N$ at the contact point. Then the contact direction at $z_0$ is $N(z_0)r$. 
\end{lemma}

{\it Proof.} See \cite{wechselberger2020}.

%

\begin{remark} We may also construct local projections onto the contact direction by selecting a nonzero row $l$ of adj$(D\!f N)$. Note that if rank$(D\!f N) < n- k - 1$, then adj$(D\!f N) = 0$, so these results are not generalizable to the case of higher-dimensional tangencies of the fast fibers with the critical manifold. \end{remark}


\section{Contact between submanifolds of $\mathbb{R}^n$} \label{sec:contact}

Our primary goal is to classify points in the contact set $F$ according to their singularity type. To do this, we must rigorously define a notion of contact between two submanifolds of $\mathbb{R}^n$. We follow the development of Izumiya et. al. \cite{izumiya2015}. 

We begin with the most elementary setting: contact between two smooth regular curves $\alpha(t)$ and $\beta(t)$ in $\mathbb{R}^2$ sharing a common point $\alpha(0) = \beta(0) = z_0$. One candidate definition of contact order is as follows:
\begin{definition} \label{def:contact2}
The curves $\alpha$ and $\beta$ make contact of order $c$ at $z_0$ if 
\begin{eqnarray*}
\alpha^{(i)}(0) &=& \beta^{(i)}(0)\text{~~for~~} i = 0, \cdots, c\\
\alpha^{(c+1)}(0) &\neq& \beta^{(c+1)}(0)  .
\end{eqnarray*}
\end{definition}

We will instead use a slightly different, though equivalent, notion of contact for two curves in the plane which turns out to be more natural for our setting. We are ultimately concerned with contact between the fast fibers and the critical manifold of system \eqref{eq1}, where the critical manifold is defined as the level set of a submersion. This motivates the following definition, where we assume that one of the curves lies inside the zero level set of a smooth function.

\begin{definition} \label{def:curves} Let $\alpha, \beta : \mathbb{R} \to \mathbb{R}^2$ define two regular curves so that the image of $\beta$ is equal to the zero set of a smooth function $F: \mathbb{R}^2 \to \mathbb{R}$ and $\alpha(0) = \beta(0)$. Note that $ (F\circ \alpha): \mathbb{R} \to \mathbb{R}$. We say that $\alpha$ and $\beta$ {\it have contact of order $c$ at $t = 0$} if 
\begin{eqnarray*}
(F\circ \alpha)^{(i)}(0) &=& 0 \text{~~for~~} i = 0, \cdots, c\\
(F\circ \alpha)^{(c+1)}(0) &\neq& 0.
\end{eqnarray*}
\end{definition}

{\it Example: contact of order one between curves.}\\
 Suppose $\alpha(t)$ and $\beta(t)$ have contact of order 1, with $F$ as in Def. \ref{def:curves}. Then
\begin{eqnarray*}
F(\alpha(0)) &=& 0\\
D\!F  (\alpha(0)) \alpha'(0) &=& 0\\
D^2\!F(\alpha'(0),\alpha'(0)) + D\!F  (\alpha(0)) \alpha''(0) &\neq & 0.
\end{eqnarray*}

The first condition specifies that the contact point $\alpha(0)$ lies in the zero set of $F$ (i.e. on the curve $\beta$). The second condition specifies that the tangent vector $\alpha'(0)$ of $\alpha$ at the contact point lies inside the tangent space of the curve $\beta$, given by ker $ D\!F  $, at the contact point.\\

 To compare the third condition, first observe that  $\beta(t)$ lies inside the zero level set of $F$ by construction, and thus $F(\beta(t)) = 0$ over an interval of $t$. Differentiating both sides twice and evaluating the result at $t = 0$, we obtain the identity
\begin{eqnarray*}
D^2\!F(\beta'(0),\beta'(0)) &=& - D\!F  (\beta(0)) \beta''(0).
\end{eqnarray*}

Using this identity, $\alpha(0) = \beta(0)$, and $\alpha'(0) = \beta'(0)$ in the third condition, we have
\begin{eqnarray*}
D\!F  (\beta(0))(\alpha''(0) - \beta''(0)) &\neq& 0,
\end{eqnarray*}
and thus
\begin{eqnarray*}
\alpha''(0) &\neq&  \beta''(0).
\end{eqnarray*}

Similar computations can be used to demonstrate the equivalence of these two definitions of contact at higher orders. \\

Our definition of contact is well-defined with respect to smooth reparametrizations of $\alpha(t)$ and $\beta(t)$. Less obvious technical issues that must be resolved are that (i) contact of a particular order should not depend on a particular choice of function $F$, and (ii) the notion of contact is inherently local, so the domains of $\alpha$ and $F$ should not matter outside of small neighborhoods of the contact point $z_0$.\\

\begin{remark}
These issues motivate the use of {\it germs of smooth functions} and {\it jet spaces}, which are the natural objects that we will use to define contact in the general case. We give a broad description of these terms, and relegate proper definitions to the Appendix (\ref{appendix}). Fix a basepoint $z_0 \in \mathbb{R}^n$. The germ of a function $f$ at  $z_0$ is defined from the equivalence class of all smooth functions $g: \mathbb{R}^n \to \mathbb{R}^m$ that are equal to $f$ on a common neighborhood of $z_0$. The collection of germs has the structure of a ring in the space of smooth functions. We may define the {\it $k$-jet space at $z_0$}, denoted $J^k(n,m)$, by taking a quotient of this ring by the ideal of all germs that vanish to order $k$. The $k$-jet space may be identified with the set of polynomials of total degree less than or equal to $k$. The $k$-jet of a germ $f$, denoted $J^k f$, is, roughly speaking, the element of $J^k(n,m)$ that may be identified with the truncated Taylor polynomial of order $k$ under some suitable local coordinate transformation. The natural advantage of using $k$-jets is that contact can be defined in a coordinate-independent manner.
\end{remark}

We now rigorously define contact between two submanifolds of $\mathbb{R}^n$, following the treatment of Montaldi \cite{montaldi1983,montaldi1986,montaldi1991} and the presentation of Izumiya et. al. \cite{izumiya2015}. The first step is to define a suitable generalization of contact order. 

\begin{definition} Suppose $(M_1,N_1)$ and $(M_2,N_2)$ are two pairs of submanifolds of $\mathbb{R}^n$ with dim$(M_i) = m$ and dim$(N_i) = d$. We say that {\it the contact of $M_1$ and $N_1$ at $y_1$ is of the same type as the contact of $M_2$ and $N_2$ at $y_2$} if there is a germ of a diffeomorphism $\Phi: (\mathbb{R}^n,y_1) \to (\mathbb{R}^n,y_2)$ so that $\Phi(M_1) = M_2$ and $\Phi(N_1) = N_2$. \end{definition}

Our objective is to relate this `generalized contact order' to a generalized version of the map $F \circ \alpha$ in Def. \ref{def:curves} in a suitable setting. We have

\begin{definition} \label{def:contactmap} Suppose that a submanifold $M \subset \mathbb{R}^n$ is given locally as the image of some immersion-germ $g:(M,x) \to (\mathbb{R}^n,0)$ and another submanifold $N \subset \mathbb{R}^n$ is given by the zero set of some submersion-germ $f: (\mathbb{R}^n,0) \to (\mathbb{R}^k,0)$. The {\it contact map} of $M$ and $N$ near $x$ is the germ of the composite map $f \circ g$ near $x$. \end{definition}

This definition should be compared with Def. \ref{def:curves}. The regular curve $\alpha$ defines an immersion in small neighborhoods of the contact point. 

\begin{remark} For contact between a $k$-dimensional submanifold $M$ and a regular (1-dimensional) curve $\alpha$, Wechselberger uses a more general version of Def. \ref{def:contact2} in \cite{wechselberger2020}. It is much less straightforward to compute the contact order using that definition, especially for contact of order larger than one. In the next section we will show how Def. \ref{def:contactmap} can be used to give a definition of contact order which is easier to compute. The equivalence of these two definitions follows from the characterization of the tangent space via equivalence classes of curves. \end{remark}

\begin{remark} For submanifolds of {\it equal} dimension, analogous definitions of contact have been used to study contact between locally defined center manifolds \cite{wan1977}, between center manifolds and center subspaces \cite{murdock2003}, and to define jet bundles \cite{olver1986}. \end{remark}

We now present two important results. 

\begin{lemma}\label{lemma:contact} For any pair of submanifolds in $\mathbb{R}^n$, the contact class of the contact map depends only on the submanifold-germs themselves and not on the choice of submersion and immersion germs (and therefore not on the contact map). \end{lemma}

{\it Proof.} See \cite{montaldi1986} or Lemmas 4.1 and 4.2 in \cite{izumiya2015}.\\

This lemma ensures that `the' contact map-germ of two submanifolds is well-defined. The contact class of a smooth germ is the equivalence class of all germs whose zero sets are diffeomorphic (the complete definition is given in the appendix). 

\begin{theorem} \label{thm:contact} 
Suppose $g_1,g_2: M_{1,2} \to \mathbb{R}^n$  are immersion-germs and $f_{1},f_2: \mathbb{R}^n \to \mathbb{R}^k$ are submersion-germs (with $N_{1,2} = f_{1,2}^{-1}(0)$). Then the pairs $(M_1,N_1)$ and $(M_2,N_2)$ have the same contact type iff $f_1 \circ g_1$ and $f_2 \circ g_2$ lie in the same contact class.
 \end{theorem}

{\it Proof.} See \cite{montaldi1986} or Theorem 4.1 in \cite{izumiya2015}.\\

This result provides the required connection to classical singularity theory: the contact class of a pair of submanifolds is completely determined by the singularities of the contact map between them.

\subsection{$A_c$ contact singularities and their unfoldings} \label{sec:ac}

We are finally in a position to consider the main setting of this paper: the contact of a curve $\alpha: \mathbb{R}\to \mathbb{R}^n$ with a submanifold given by the zero set of a submersion $f: \mathbb{R}^n \to \mathbb{R}^{n-k}$ (with $1 \leq k < n$). The contact map is $ f\circ\alpha: \mathbb{R} \to \mathbb{R}^{n-k}$. \\

By Theorem \ref{thm:contact}, the contact type is well-defined by the contact class of $f\circ{\alpha}$. We are interested in {\it stable} maps with respect to the contact class, i.e equivalence classes of maps which have versal unfoldings of finite codimension (see \cite{izumiya2015} Section 3.8 for complete definitions of stable maps and versal unfoldings and Theorem 3.9 for the relationship between these two notions). The maps $h: \mathbb{R} \to \mathbb{R}^{d \geq 1}$ which are stable with respect to the contact class are the well-understood and classified $A_c$ singularities (see for eg. \cite{guseinzade, izumiya2015}). 

\begin{definition}
(See Sec. 11.1 of \cite{guseinzade}) A critical point $p$ of a smooth function $f: \mathbb{R}^n \to \mathbb{R}$ is {\it of type $A_c$} if $f$ is locally equivalent near $p$ to $x_1^{c+1} + x_2^2 + \cdots +x_n^2$.
\end{definition}

Stable maps are finitely-determined \cite{mather1969b}, so each stable map is contact-equivalent (see Def. \ref{def:contactequiv} in the Appendix) to the germ of a map $\tilde{h}(t) = (t^{c+1},0,\cdots, 0)$ for some constant $c$. In this way we readily obtain a suitable analogue of the $A_c$ singularity classes for contact equivalence:

\begin{definition} 
A smooth map $h: \mathbb{R} \to \mathbb{R}^{n-k}$ {\it has an $A_c$ singularity} if  it is contact-equivalent to $\tilde{h}(t) = (t^{c+1},0,\cdots, 0)$. The derivative conditions for an $A_c$ singularity are
\begin{eqnarray*}
h^{(i)}(0) &=& 0 \text{ for } i = 0, \cdots, c\\
h^{(c+1)} &\neq& 0.
\end{eqnarray*}
\end{definition}

\begin{definition}\label{def:akcontact}
If $h = f \circ \alpha$ is a contact map between a $k$-dimensional submanifold given by the zero set of a submersion $f: \mathbb{R}^n \to \mathbb{R}^{n-k}$ and a curve  $\alpha: \mathbb{R}\to \mathbb{R}^n$, we say that the submanifold and the curve {\it make contact of order $c$ at $z_0$} if $h$ admits an $A_{c}$-singularity at $z_0$.
\end{definition}

%
%
%
%
%

We demonstrate the computability of Def. \ref{def:akcontact} with an example. Let $\alpha(t): \mathbb{R}\to \mathbb{R}^n$ be a curve in $\mathbb{R}^n$ and let $M$ be a $k$-dimensional submanifold in $\mathbb{R}^n$ given by the zero level set of a submersion: $M = f^{-1}(0)$. Let $\alpha(0) = z_0 \in M$ denote a contact point between $\alpha$ and $M$.

By Def. \ref{def:akcontact},  $\alpha$ makes contact of order 2 with $M$ at $z_0$ if 
\begin{eqnarray*}
(f \circ \alpha)(0) &=& 0\\
(f\circ\alpha)'(0) &=& 0\\
(f\circ\alpha)''(0) &=& 0\\
(f\circ\alpha)'''(0) &\neq& 0.	
\end{eqnarray*}

These derivatives may be evaluated using repeated applications of the chain rule:
\begin{eqnarray*}
D\!f  (z_0) \alpha'(0) &=& 0\\
D^2\!f(\alpha'(0),\alpha'(0)) + D\!f  \alpha''(0) &=& 0\\
D^3\!f(\alpha'(0),\alpha'(0),\alpha'(0)) +3 D^2\!f(\alpha'(0),\alpha''(0))+ D\!f   \alpha'''(0) &\neq& 0. \label{eq:thirdder}
\end{eqnarray*}

Note the geometric content of the first condition $D\!f  (z_0) \alpha'(0) = 0$: the tangent vector of $\alpha$  lies precisely inside the tangent space of $M$ at the contact point.

\begin{remark} {\it Multilinear maps.}  We remind the reader of the standard notation and evaluation of multilinear maps $L: \mathbb{R}^{n_1}\times \cdots \times \mathbb{R}^{n_k} \to \mathbb{R}^{n_d}$:
\begin{eqnarray*}
L(v_1, \cdots, v_{k})_{j} &=& \sum_{l_{k} = 1}^{n_k} \cdots  \sum_{l_{1} = 1}^{n_1} L_{l_1,\cdots, l_k,j} v_{1,l_1} \cdots v_{k,l_k} \qquad \mbox{for } j = 1,\cdots,n_d.
\end{eqnarray*}

For example, 
\begin{eqnarray*}
(D^2\!f(\alpha'(0),\alpha'(0)))_i &=&  \sum_{j,k=1}^{n} \frac{\partial^2 f_i}{\partial x_j \partial x_k} \alpha_j'(0)\alpha_k'(0) \qquad \mbox{for }i = 1, \cdots, n-k.
\end{eqnarray*}
\end{remark}

\subsection{Computing the contact order} \label{sec:straighten}

At points in $z_0 \in F$, the center manifold theorem provides local families of one-dimensional center manifolds $W^c(z_0)$ of the layer problem all tangent to the contact direction $Nr$ at the basepoint $z_0$. 
We may thus define regular curves passing through $z_0$ with nonzero speed. Our goal is to evaluate the derivatives of a given curve segment $\alpha(t)$ at $\alpha(0) = z_0$, in terms of derivatives of $N$ and $f$. \\

{\it The subcase $\dim S = n-1$.}\\

The factorization $h_0(z) =  N(z)f(z)$ consists of the term $N(z)$ of size $n \times 1$ and the scalar function $f(z)$. Let $z_0 \in F$ and let $B_{z_0}$ denote an open ball centered at $z_0$ such that $N(z)$ is nonzero for all $z \in B_{z_0}$ and $f(z) \neq 0$ for all $z \in B_{z_0}-S$. For points $z \in B_{z_0} - S$, the vector field is a nonzero multiple of the vector field $N(z)$. Solutions of the desingularized layer problem
\begin{eqnarray}
z' &=& N(z) \label{eq:desing2}
\end{eqnarray}

defined in $B_c$ consist of regular curves. In particular, each such curve crosses $S$ with nonzero speed. Away from $S$, the tangent vectors of the solution curves are aligned with the original vector field $h_0$ everywhere (possibly with a change in orientation).\\

In the codimension-1 case it is straightforward to compute high-order derivatives of solution curves of \eqref{eq:desing2}. Let $\alpha(t)$ denote a solution curve of \eqref{eq:desing2} with the property that $\alpha(0) = z_0 \in S$. Then

\begin{eqnarray*}
\alpha'(0) &=& N(z_0)\\
\alpha''(0) &=& D\!N(z_0) N(z_0)\\
\alpha'''(0) &=& D^2N(z_0)(N(z_0),N(z_0)) + D\!N(z_0)D\!N(z_0)N(z_0) ,
\end{eqnarray*}

etc. For example, if $z_0 \in F$, then $\alpha'(0)  = N(z_0) \neq 0$, whereas $(f \circ \alpha)'(0) = D\!f (z_0)  N(z_0) = 0$, giving contact order of at least one---as expected. \\

{\it The general case $1 \leq \dim S \leq n-1$.}\\

The technique of desingularising the layer problem does not generalise to the case $\dim S < n-1$ since $f$ is no longer a scalar function--- the components $f_i(\alpha(t))$ grow and shrink independently along a curve $\alpha(t)$ lying in the center manifold. We sidestep this technical issue by first making a local coordinate transformation which straightens the fibers locally.

\begin{lemma} \label{lemma:straighten}

In the following, we evaluate all quantities at $z_0 \in F$. Recall that $r$ is a right nullvector of $D\!f  (z_0)N(z_0)$.The defining equations for contact of order-one between the fast fiber bundle and the critical manifold are
\begin{eqnarray*}
D\!f   Nr &=& 0\\
D^2\!f(Nr,Nr) + D\!f  D\!N(Nr,r) &\neq&  0.
\end{eqnarray*}

The defining equations for contact of order-two between the fast fiber bundle and the critical manifold are
\begin{eqnarray*}
D\!f   Nr &=& 0\\
D^2\!f(Nr,Nr) + D\!f  D\!N(Nr,r) &=&  0\\
D^3\!f(Nr,Nr,Nr) +3 D^2\!f(Nr,D\!N(Nr,r))+ &&\\
 D\!f   (D^2N(Nr,Nr,r) + D\!N(D\!N(Nr,r),r))  &\neq& 0. 
\end{eqnarray*}

Generally, the defining equations for contact of order-$c$ at a point $z_0 \in F$ are
\begin{eqnarray*}
\nabla_{Nr}^{j}f|_{z_0} &=& 0   \qquad \mbox{for }j = 0,\cdots, c\\
\nabla_{Nr}^{c+1} f|_{z_0} &\neq& 0,
\end{eqnarray*}

where $\nabla_{Nr}$ denotes the directional derivative in the direction $Nr$. 
\end{lemma}

{\it Proof.} There exists a local coordinate transformation $u = L(z)$ which places the system \eqref{eq1} in standard form by straightening the fibers, with local coordinates $(u,y)$ (see \cite{fenichel1979}, \cite{wechselberger2020} for the full treatment, or \cite{jones1995} for the straightening step beginning with the system in standard form). If we let $x = M(u,y)$ denote the local inverse of $L$, we have

\begin{eqnarray*}
\begin{pmatrix} u'\\ y' \end{pmatrix} &=& \begin{pmatrix} \mathbb{O}_{k,n-k} \\ N^y(M(u,y),y) \end{pmatrix} f(M(u,y),y)\\
&=&  \begin{pmatrix} \mathbb{O}_{k,n-k} \\ I_{n-k,n-k} \end{pmatrix} \tilde{f}(u,y).
\end{eqnarray*}

Here, $u \in \mathbb{R}^k$ is the local slow variable and $y \in \mathbb{R}^{n-k}$ is the local fast variable, so that $z_0$ is placed at $(u,y) = (0,0)$. The direction of contact can be made explicit through another change of variable

\begin{eqnarray*}
y &=& N^y r v + N^y P w,
\end{eqnarray*}
 
 where $P$ is an $(n-k) \times (n-k-1)$ matrix chosen so that Col$\begin{pmatrix} r & P \end{pmatrix}$ provides a basis of $\mathbb{R}^k$. Let $l$ and $Q$ satisfy the identity
 
 \begin{eqnarray*}
\begin{pmatrix} l \\ Q \end{pmatrix} \begin{pmatrix} r & P \end{pmatrix} &=& I_{n-k,n-k}.
\end{eqnarray*}

Let $\tilde{r} = N^y r$ and $\tilde{P} = N^y P$. After these two coordinate transformations, we have

\begin{eqnarray*}
\begin{pmatrix} u' \\ v' \\ w' \end{pmatrix} &=& \begin{pmatrix} \mathbb{O}_{k,n-k} \\ l f(M(u,\tilde{r}v+\tilde{P}w),\tilde{r}v+\tilde{P}w) \\ Q f(M(u,\tilde{r}v+\tilde{P}w),\tilde{r}v+\tilde{P}w)\end{pmatrix} \\
&=& \begin{pmatrix} \mathbb{O}_{k,n-k} \\ I_{n-k,n-k} \end{pmatrix} \begin{pmatrix}  l f(M(u,\tilde{r}v+\tilde{P}w),\tilde{r}v+\tilde{P}w) \\ Q f(M(u,\tilde{r}v+\tilde{P}w),\tilde{r}v+\tilde{P}w)\end{pmatrix}\\
&=& \begin{pmatrix} \mathbb{O}_{k,n-k} \\ I_{n-k,n-k} \end{pmatrix} \tilde{f}(u,v,w),
\end{eqnarray*}

where by slight abuse of notation we use $\tilde{f}$ again to denote the second factor. We observe that in particular the one-dimensional center manifold has been locally straightened, containing the regular curve

\begin{eqnarray*}
\tilde{\alpha}(t) &=& \begin{pmatrix} \mathbb{O}_{k,1} \\ 1 \\ \mathbb{O}_{n-k-1,1}  \end{pmatrix} t.
\end{eqnarray*}

The straightening transformations simultaneously deform the critical manifold locally, reflected in the transformation from $f$ to $\tilde{f}$. The key point is that this sequence of coordinate transformations {\it preserves the contact order} (by Theorem \ref{thm:contact}).\\

It remains to compute the defining equations for the $A_c$  singularity classes of the deformed contact map $\tilde{f} \circ \tilde{\alpha}$. Observe that
\begin{eqnarray*}
(\tilde{f} \circ \tilde{\alpha})(0) &=& 0\\
(\tilde{f} \circ \tilde{\alpha})^{(k)}(0) &=& D^{(k)}\tilde{f}(0) (\tilde{\alpha}'(0), \cdots, \tilde{\alpha}'(0)) \\
&=& D_{\underbrace{v\cdots v}_{k \text{ times}} }\tilde{f}(0)
\end{eqnarray*}

since $\tilde{\alpha}^{(k)}(0) = 0$ for $k \geq 2$, and the subscript in the last line denotes that the partial derivative is taken $k$ times with respect to the center variable $v$. We compute the first two derivatives of the contact map. We have

\begin{eqnarray*}
D_v \tilde{f}(M(u,\tilde{r}v + \tilde{P}w),\tilde{r}v + \tilde{P}w) &=& \begin{pmatrix} l \\ Q \end{pmatrix}D_xf D_yM \tilde{r} + D_y f \tilde{r} \\
&=&  \begin{pmatrix} l \\ Q \end{pmatrix}(D_x f N_x (N^y)^{-1} + D_yf) \tilde{r}\\
&=& \begin{pmatrix} l \\ Q \end{pmatrix} (D_x f N_x + D_yf N^y )(N^y)^{-1} N^y r\\
&=& \begin{pmatrix} l \\ Q \end{pmatrix}D\!f   N r.
\end{eqnarray*}

Generally, for any differentiable function 
\begin{eqnarray*}
\zeta &=& \zeta(M(u,\tilde{r}v + \tilde{P}w),\tilde{r}v + \tilde{P}w), 
\end{eqnarray*}
we have
\begin{eqnarray*}
D_v \zeta &=& D\zeta Nr.
\end{eqnarray*}

This rule can be used to generate derivatives of arbitrary order. We have for instance the second-order derivative
\begin{eqnarray*}
D_{vv} \tilde{f}(M(u,\tilde{r}v + \tilde{P}w),\tilde{r}v + \tilde{P}w) &=& \begin{pmatrix} l \\ Q \end{pmatrix}(D^2\!f (Nr,Nr) + D\!f   D\!N(Nr,r)),
\end{eqnarray*}
the third-order derivative
\begin{eqnarray*}
D_{vvv} \tilde{f}(M(u,\tilde{r}v + \tilde{P}w),\tilde{r}v + \tilde{P}w) &=& \begin{pmatrix} l \\ Q \end{pmatrix}(D^3\!f(Nr,Nr,Nr) +3 D^2\!f(Nr,D\!N(Nr,r))+ \\
&& D\!f   (D^2N(Nr,Nr,r) + D\!N(D\!N(Nr,r),r))),
\end{eqnarray*}

and so on. The matrix $\begin{pmatrix} l \\ Q \end{pmatrix}$ stores the dual basis, and therefore does not affect the corresponding (non)zero conditions in the associated defining equations. \\

To prove the final statement of the lemma, recall the definition of the directional derivative: for a test function $\xi$, we have
\begin{eqnarray*}
\nabla_{Nr}\xi &=& D\xi (Nr).
\end{eqnarray*}
The proof follows by repeated applications of the chain rule.
 \hfill$\Box$\\

We emphasize that a straightening transformation is {\it not} required to check the conditions of Lemma \ref{lemma:straighten}; it is only used as an ingredient in the proof. The straightening transformation in the lemma above has been used to compute an unfolding of contact points of order-one, with the slow variables unfolding the contact point \cite{wechselberger2020}. The promised identification of the nondegeneracy condition for the $A_c$ singularity of the contact map, and the corresponding codimension-$c$ folded singularity nondegeneracy condition, is now established, up to a trivial projection along the center flow. \\

\subsection{Slow unfoldings of contact singularities; contact folds} \label{sec:folds}

From the point of view of geometric singular perturbation theory, the dynamical relevance of loss of normal hyperbolicity of the critical manifold is only manifested when the unfoldings occur {\it under local variation of the slow variables}. For example, for the case of standard slow-fast systems \eqref{eq:standard} satisfying the fold conditions \eqref{eq:foldconds} plus an additional `slow dynamics transversality' condition $g(0,0,0) \neq 0$, there exists \cite{krupa2001} a smooth, local coordinate change $\phi(x,y) = (\xi,\eta)$ in which \eqref{eq:standard} is given by
\begin{eqnarray*}
\frac{d \xi}{dt} &=& \eta + \xi^2 + \mathcal{O}(\xi^2, \xi \eta, \eta^2,\eps)\\
\frac{d \eta}{dt}&=& \eps( \pm 1 + \mathcal{O}(\xi,\eta,\eps)).
\end{eqnarray*}

Slow unfoldings of folded nodes admit scenarios where the slow flow may cross fold points transversely. This is the basic ingredient in constructing persistent nontrivial connections between attracting and repelling slow manifolds. 

%

\begin{remark}
The analysis of normal forms of folded nodes has been extended to higher dimensions. In the $n =3$ case of two slow variables and one fast variable, a suitable time-rescaling of the layer problem can be recast as a Riccati equation to leading order \cite{szmolyan2004}. Remarkably, this observation holds for the $k$-slow $m$-fast case (with $k \geq 2$ and $m \geq 1$) as well, allowing the theory developed in $\mathbb{R}^3$ to be extended to arbitrary dimensions \cite{wechselberger2012}. 
\end{remark}

 We now consider computable criteria for such restricted unfoldings in the more general case of contact singularities.  
 
 \begin{definition} Assume $\eqref{eq1}$ admits a contact point $z_0 \in F$. We say that $z_0$ is a {\it contact fold} if it is a contact point of order one that admits a versal (codimension-one) unfolding under local variation of the slow variables. 
\end{definition}

The question is how to write down coordinate-independent defining equations for contact folds {\it without} having to compute explicit local coordinate changes everywhere along the critical manifold. The first main result of our paper is proven by following a similar procedure   to that of Lemma \ref{lemma:straighten}: we apply a local coordinate transformation of \eqref{eq1} in which local slow variables may be identified and extract a geometric condition from the transformed system. 

\begin{lemma}\label{lemma:fold}
Defining equations for a contact fold at $z_0 \in F$ are:
\begin{itemize}
\item[(a)] $\text{rank}(D\!f  (z_0)) = n-k$
\item[(b)] $l(D^2\!f(Nr,Nr) + D\!f   D\!N(Nr,r))|_{z_0} \neq 0,$
\end{itemize}

where $l,r$ are nonzero left and right nullvectors of $D\!f  (z_0)N(z_0)$. \end{lemma}

{\it Proof.} The first local regularity condition is satisfied immediately when $f$ is a submersion. We follow the proof given in \cite{wechselberger2020} to compute the defining equation.  From the proof of Lemma \ref{lemma:straighten} we recall the transformed vector field with locally straightened fibers
\begin{eqnarray}
\begin{pmatrix} u' \\ v'\\ w' \end{pmatrix} &=&  \begin{pmatrix} \mathbb{O}_{k,1} \\ \tilde{l} \\ \tilde{Q} \end{pmatrix} \tilde{f}(u,\tilde{r} v + \tilde{P} w) \nonumber\\
&=& \begin{pmatrix} \mathbb{O}_{k,1} \\ l \\ Q \end{pmatrix} f(M(u,\tilde{r} v + \tilde{P}w),\tilde{r} v + \tilde{P}w), \label{eq:straightened}
\end{eqnarray}

where $(u,v,w) \in \mathbb{R}^k \times \mathbb{R} \times \mathbb{R}^{n-k-1}$. \\

The Jacobian $\tilde{J} = D\tilde{h}_0$ along $S$ is given by
\begin{eqnarray*}
\tilde{J}|_S &=& \begin{pmatrix} \mathbb{O}_{k,k} & \mathbb{O}_{k,1} & \mathbb{O}_{k,n-k-1} \\ 
l D_xf (D_x L)^{-1} & l D\!f   N r & l D\!f   NP\\ 
Q D_x f (D_xL)^{-1} & QD\!f   N r & QD\!f   NP \end{pmatrix}.
\end{eqnarray*}

On $F$, the (2,2), (3,2), and $(2,3)$ (block) entries are further annihilated because $l$ and $r$ are precisely the nullvectors of $D\!f   N$ on the set $F$ of contact points of $S$:
\begin{eqnarray*}
\tilde{J}|_F &=& \begin{pmatrix} \mathbb{O}_{k,k} & \mathbb{O}_{k,1} & \mathbb{O}_{k,n-k-1} \\ 
l D_xf (D_x L)^{-1} & \mathbb{O}_{1,1} & \mathbb{O}_{1,n-k-1} \\ 
Q D_x f (D_xL)^{-1} & \mathbb{O}_{n-k-1,1} & QD\!f   NP \end{pmatrix}.
\end{eqnarray*}

 Near $F$ we expand the right-hand side of the (one-dimensional) $v'$ equation. We have
\begin{eqnarray*}
v' &=& \tilde{l} \tilde{f}(u,\tilde{r}v + \tilde{P}w)\\
&=& lf(M(u,\tilde{r}v+\tilde{P}w),\tilde{r}v+\tilde{P}w)\\
&=&  l D_x f (D_x L)^{-1} u + l (D^2\!f (Nr,Nr) + D\!f   D\!N(Nr,r))v^2 + \cdots ,
\end{eqnarray*} 

where we ignore the remaining cross-terms of order two and the higher-order terms.\\

The coefficient of the vector-valued component $u$ is $l D_x f (D_x L)^{-1}$ which is nontrivial since rank $D\!f  (z_0)  = n-k$. Thus $l D_x f (D_x L)^{-1} u$ plays the role of an unfolding parameter, but with the parameter axis lying along a nullvector of  $l D_x f (D_x L)^{-1}$.  \hfill $\Box$

\subsection{Contact cusps} \label{sec:cusps}
\begin{definition} Assume $\eqref{eq1}$ exhibits a contact point $z_0 \in F$. We say that $z_0$ is a {\it contact cusp} if it is a contact point of order two that admits a versal (codimension-two) unfolding under local variation of the slow variables. 
\end{definition}

We recall the generic criteria for the unfolding of a cusp point \cite{kuznetsov}. Consider the smooth vector field
\begin{eqnarray}
\dot{x} &=& f(x,\alpha_1,\alpha_2), \label{eq:cusp0}
\end{eqnarray}

$x,\alpha_1,\alpha_2 \in \mathbb{R}$, with an isolated equilibrium point at $p = (x,\alpha_1,\alpha_2) = (0,0,0)$. Assume the following:
\begin{itemize}
\item $D_xf(p) = D_{xx}(p) = 0$ 
\item Nondegeneracy condition: $D_{xxx}f(p) \neq 0 $
\item Parameter transversality condition: $\det \begin{pmatrix} D_{\alpha_1}f(p) & D_{\alpha_2}f(p) \\
D_{x\alpha_1}f(p) & D_{x\alpha_2}f(p) \end{pmatrix} \neq 0$.
\end{itemize}

Then we can find smooth invertible coordinate change in the extended phase space so that the system \eqref{eq:cusp0} is transformed into
\begin{eqnarray*}
\dot{\eta} &=& \beta_1 + \beta_2 \eta + \eta^3 + \mathscr{O}(\eta^4).
\end{eqnarray*}

These nondegeneracy and parameter transversality conditions can be expressed more compactly by specifying instead that the map $F: \mathbb{R}^3 \to \mathbb{R}^3$ defined by
\begin{eqnarray*}
F: (x,\alpha_1,\alpha_2) \mapsto (f,D_xf,D_{xx}f)(x,\alpha_1,\alpha_2)
\end{eqnarray*}

be regular at the cusp point. 

\begin{remark}
In two-parameter families of $n$-dimensional flows $\dot{x} = f(x,\alpha_1,\alpha_2)$ with $x \in \mathbb{R}^n$, the corresponding conditions for a cusp bifurcation are that $D_x\!f$ has one simple zero eigenvalue and $n-1$ eigenvalues with nonzero real part. Then there exist coordinate transformations locally placing the vector field in the normal form
\begin{eqnarray*}
\dot{u} &=& \beta_1 + \beta_2 u \pm u^3\\
\dot{y} &=& Ay,
\end{eqnarray*}

where $\beta_1,\beta_2$ are parameters, $(u,y) \in \mathbb{R}\times \mathbb{R}^{n-1}$ and $A$ is an $(n-1)\times (n-1)$ hyperbolic matrix \cite{kuznetsov}.
\end{remark}

We now state and prove the analogous result to Lemma \ref{lemma:fold} for contact cusps of nonstandard slow-fast systems.

\begin{lemma} \label{lemma:cusp}
The defining equation and genericity conditions for a contact cusp at $z_0 \in F$ are:
\begin{itemize}
\item[(a)] $\text{rank}(D\!f  (z_0)) = n-k$
\item[(b)] $l(D^2\!f(Nr,Nr) + D\!f   D\!N(Nr,r))|_{z_0} =0$
\item[(c)] $l\cdot ( D^3\!f(Nr,Nr,Nr) + 3D^2\!f (D\!N(Nr,r),Nr) +$  \\ $D\!f   D^2N(Nr,Nr,r) + D\!f   D\!N(D\!N(Nr,r),r))|_{z_0} \neq  0$
\item[(d)]  The $2\times n$ matrix 
\begin{eqnarray*}
C_0(z_0) &=&\left.\begin{pmatrix} l D\!f   \\ l (D^2\!f(Nr,I) +  D\!f   D\!N(r,I))\end{pmatrix}\right\vert_{z_0}
\end{eqnarray*}
has full rank of $2$.
\end{itemize}

Here, $l$ and $r$ denote nonzero left and right nullvectors of $D\!f  (z_0)N(z_0)$.
\end{lemma}

\begin{remark}
The multilinear maps in the final item of Lemma \ref{lemma:cusp} are defined columnwise. For example,
\begin{eqnarray*}
(D^2\!f(Nr,I))_j &=&  D^2\!f(Nr,e_j)  \qquad \mbox{for }j = 1,\cdots, n,
\end{eqnarray*}
where $e_j$ is the $j$th unit vector. The term $D\!fD\!N(r,I)$ is also defined columnwise. The term $(D^2\!f(Nr,I) +  D\!f   D\!N(r,I))$ is therefore of size $(n-k) \times n$.

\end{remark}

{\it Proof:}
The first local regularity condition is satisfied immediately when $f$ is a submersion. Taylor-expanding the right-hand side of the $v'$ equation (see Eq. \eqref{eq:straightened}) near $(u,v,w) = (0,0,0)$, we have

\begin{eqnarray*}
(lf)(u,v,w) &=& (lf)(0,0,0) + D\!(lf) z  + D^2\!(lf)(z,z) + \cdots\\
&=& 0 + D_u(lf) u + D_v (lf) v + D_w (lf) w + (D_{uv}(lf) u+D_{wv}(lf) w)v + \\
&& \frac{1}{2}D_{vv}(lf) v^2 + \frac{1}{6}D_{vvv}(lf) v^3 + \cdots.
\end{eqnarray*}

In our setting, we have $D_v (lf)= 0$ and $D_u (lf)\neq 0$ on the set $F$.  Using this expansion, we can read off the defining equations for a cusp at a point $z_0 \in F$:

\begin{itemize}
\item[(i)] $D_{vv}(lf)(z_0) = 0$. 
\item[(ii)] $D_{vvv}(lf)(z_0) \neq 0$.
\item[(iii)] The $2\times k$ matrix
\begin{eqnarray*}
C(z) &=& \left.\begin{pmatrix} D_u(lf)  \\
D_{uv}(lf)  \end{pmatrix}\right\vert_z
\end{eqnarray*}
has full rank at the contact point: rank $C(z_0) = 2$.
\end{itemize}

The third condition provides two unfolding directions lying along the critical manifold. Note that we require $k \geq 2$. For contact between the fast fiber bundle and the critical manifold to be defined, we therefore require the system to be at least three-dimensional with a two-dimensional critical manifold.\\

These three conditions should be compared to the standard defining equations of the standard generic cusp. In particular,  the tangency direction $v$ plays the role of the unfolding variable $x$, and two linearly independent combinations of the remaining slow variables $u_1,u_2$ play the role of the unfolding parameters $\alpha_1,\alpha_2$.\\

{\it Evaluating the nondegeneracy condition.} \\

Differentiate three times and use the chain rule:
\begin{eqnarray*}
l D_{vvv}(f)|_{z_0} &=& l D_{vv} (D\!f   Nr)|_{z_0}\\
 &=& l D_{v} (D^2\!f(Nr,Nr)+ D\!f   D\!N(Nr,r))|_{z_0} \\
&=& l\cdot ( D^3\!f(Nr,Nr,Nr) + 3D^2\!f (D\!N(Nr,r),Nr) +  \\
&& + D\!f   (D^2N(Nr,Nr,r) + D\!N(D\!N(Nr,r),r)))|_{z_0}.
\end{eqnarray*}

{\it Evaluating the transversality condition.}\\

We have $l D_u f = l D_xf (D_xL)^{-1} = l D\!f   \hat{I} (D_xL)^{-1}$, where $\hat{I} =  \begin{pmatrix} I_{k,k} \\ O_{n-k,k} \end{pmatrix}$. Then

We have
\begin{eqnarray*}
\text{rank}(C(z_0)) &=& \text{rank}\left.\begin{pmatrix} l D_u f \\ l D_{uv} f \end{pmatrix}\right\vert_{z_0}\\
&=&  \text{rank}\left.\begin{pmatrix} l D\!f  \hat{I} (D_x L)^{-1} \\ l D_u(D\!f   Nr) \end{pmatrix}\right\vert_{z_0}\\
&=&  \text{rank}\left.\begin{pmatrix} l D\!f  \hat{I} (D_x L)^{-1} \\ l (D^2\!f(Nr,\hat{I}(D_x L)^{-1}) +  D\!f   D\!N(r,\hat{I}(D_x L)^{-1}))\end{pmatrix}\right\vert_{z_0}\\
&=&  \text{rank}\left.\begin{pmatrix} l D\!f   \\ l (D^2\!f(Nr,I) +  D\!f   D\!N(r,I))\end{pmatrix}\right\vert_{z_0}\\
&=& \text{rank} ~C_0(z_0),
\end{eqnarray*}

where the penultimate line follows from right-factoring the $n\times k$ full-rank matrix $\hat{I}(D_x L)^{-1}$ from both block rows of the original matrix. \hfill $\Box$\\

\subsection{Slow unfoldings of contact points of arbitrary order}
\eIL{
In this paper we focus on biochemical examples having contact folds and contact cusps, but we can extend the arguments in the proof of Lemma \ref{lemma:cusp} to derive defining equations for generic unfoldings of contact singularities of arbitrary order.

\begin{definition} Assume $\eqref{eq1}$ exhibits a contact point $z_0 \in F$ of order $c$. The contact point $z_0$ is {\it slow-generic} if it admits a versal (codimension-$c$) unfolding under local variation of the slow variables. 
\end{definition}

\begin{lemma} \label{lemma:arbitrary}
The defining equations for a slow-generic contact singularity  at $z_0 \in F$ of order $c$ are:
\begin{itemize}
\item[(a)] $\text{rank}(D\!f  (z_0)) = n-k$
\item[(b)] \begin{eqnarray*}
l \, \nabla_{Nr}^{j}f|_{z_0} &=& 0   \qquad \mbox{for }j = 0,\cdots, c\\
l \, \nabla_{Nr}^{c+1} f|_{z_0} &\neq& 0,
\end{eqnarray*}
\item[(c)]  The $c\times n$ matrix 
\begin{eqnarray*}
C_0(z_0) &=&\left.\begin{pmatrix} l D(\nabla_{Nr}^{0} f)   \\ \vdots  \\ l  D(\nabla_{Nr}^{c-1}f) \end{pmatrix}\right\vert_{z_0}
\end{eqnarray*}
has rank equal to $c$.
\end{itemize}

Here, $l$ and $r$ denote nonzero left and right nullvectors of $D\!f  (z_0)N(z_0)$.
\end{lemma}

{\it Proof:} Part (a) follows immediately when $f$ is a submersion. Part (b) follows from the identical Taylor series expansion as that given in Lemma \ref{lemma:cusp}, where the derivative conditions have been translated into directional derivatives using the final part of Lemma  \ref{lemma:straighten}. Finally, we compare the Taylor series to a versal unfolding 
\begin{eqnarray*}
f(x,\alpha_1,\cdots,\alpha_c) &=& x^{c+1} + \alpha_1 x^{c-1} + \cdots + \alpha_{c-1} x + \alpha_c
\end{eqnarray*}

of the generic $A_c$ singularity (see for eg. Part I Chapter 3 of \cite{guseinzade2}). The corresponding slow unfolding parameters are chosen from linear combinations of the coefficients of the terms $uv^{d}$ in the Taylor expansion (for $0 \leq d \leq c-1$). Therefore, we require that the $c\times k$ matrix
\begin{eqnarray*}
C(z) &=& \left. \begin{pmatrix} lD_uf \\ \vdots \\ l D_uD_{ \underbrace{v\cdots v}_{(c-1) \text{ times}}} f \end{pmatrix} \right\vert_z
\end{eqnarray*}
have rank $c$ at the contact point $z_0$. Note that this implies that $k \geq c$, i.e. we require at least $c$ slow variables for a generic slow unfolding of a contact singularity of order $c$. The formula in part (c) follows by the identical right-factorisation argument used in Lemma \ref{lemma:cusp}. \hfill $\Box$}

\section{Examples} \label{sec:examples}

We now use Lemmas \ref{lemma:fold} and \ref{lemma:cusp} in a series of three-dimensional multiple-timescale systems having a two-dimensional critical manifold, as described in the introduction. 

\subsection{Standard slow-fast systems.}\label{sec:standardcusp}

\subsubsection{The cusp normal form.}
Consider the normal form of the singularly perturbed cusp in  $\mathbb{R}^3$ in the standard case \cite{broer2013}: 

\begin{eqnarray*}
x' &=& \eps (1 + \bigO(x,y,z,\eps))\\
y' &=& \eps \bigO(x,y,z,\eps)\\
z' &=& (z^3+yz+x) + \bigO(\eps,xz,z^4).
\end{eqnarray*}

Here the slow variables are $x,y$ and the fast variable is $z$. In terms of the $Nf$-splitting we have the right-hand side of the layer problem given by

\begin{eqnarray*}
h &=& Nf = \begin{pmatrix} 0 \\ 0 \\ 1 \end{pmatrix} (x+yz+z^3).
\end{eqnarray*}

We check the conditions for a contact cusp. We first check Lemma \ref{lemma:cusp}(b). Note that 
\begin{eqnarray*}
D\!f   N &=& y + 3z^2.
\end{eqnarray*}

The parabola $F = \{(x,y,z): y+3z^2 = 0\} \cap S$ on the cusp surface $S = \{(x,y,z): x + yz + z^3 = 0\}$ consists of contact points of at least order one. In particular we have that $l = r = 1$ for this problem and thus

\begin{eqnarray*}
l (D^2\!f(Nr,Nr) + D\!f   D\!N(Nr,r)) &=&N^T D^2\!f N + D\!f   D\!N N\\
&=& \begin{pmatrix} 0 & 0 & 1 \end{pmatrix}  \begin{pmatrix} 0 & 0 & 0 \\ 0 & 0 & 1 \\ 0 & 1 & 6z \end{pmatrix} \begin{pmatrix} 0 \\ 0 \\ 1 \end{pmatrix} + 0\\
&=& 6z,
\end{eqnarray*}

whence every point on $F$ except the point $z_0 = (0,0,0)$ is a fold point.\\

The condition in  Lemma \ref{lemma:cusp}(c)  reduces to computing

\begin{eqnarray*}
l D^3\!f (N,N,N) &=& \sum_{j,k,l = 1}^3 \frac{\partial^3 f}{\partial z_j \partial z_k \partial z_m} N_j N_k N_m\\
&=& 6 \neq 0.
\end{eqnarray*}

Finally, we check the transversality condition Lemma \ref{lemma:cusp}(d). We have

\begin{eqnarray*}
C_0(z_0) &=& \left.\begin{pmatrix} D\!f   \\ D^2\!f N + D\!f   D\!N  \end{pmatrix} \right\vert_{z_0}\\
&=& \begin{pmatrix} 1 & 0 & 0 \\  0 & 1 & 0\end{pmatrix},
\end{eqnarray*}

as expected.\\

Defining equations for a contact cusp in standard slow-fast systems are given by a trivial subcase of Lemma \ref{lemma:cusp}. Consider the layer problem

\begin{eqnarray*}
\begin{pmatrix} x' \\ y' \\ z' \end{pmatrix} &=& \begin{pmatrix} 0 \\ 0 \\ 1 \end{pmatrix} g(x,y,z),
\end{eqnarray*}

which is given in standard form. Then $N(x,y,z) = \begin{pmatrix} 0 & 0 & 1 \end{pmatrix}^T$. Then at a test point $p_0 = (x_0,y_0,z_0)$ the (non)degeneracy conditions from \ref{lemma:cusp} (b) and (c) on the derivatives become

\begin{eqnarray*}
g = D_z g = D_{zz}g &=& 0\\
D_{zzz}g &\neq& 0,
\end{eqnarray*}

whereas the test matrix for the transversality condition (Lemma \ref{lemma:cusp}(d)) is 

\begin{eqnarray*}
C_0(z_0) &=& \begin{pmatrix} D\!f   \\ D^2\!f N + D\!f  D\!N \end{pmatrix}\\
&=& \begin{pmatrix} D_x g & D_y g & D_z g \\ D_{xz} g& D_{yz}g & D_{zz} g\end{pmatrix}\\
&=& \begin{pmatrix} D_x g & D_y g & 0 \\ D_{xz} g& D_{yz}g & 0 \end{pmatrix},
\end{eqnarray*}

giving the transversality condition
\begin{eqnarray*}
\det \begin{pmatrix}  D_x g& D_y g\\ D_{xz} g& D_{yz}g \end{pmatrix}&=&\begin{pmatrix} D_x g \\ D_yg \end{pmatrix} \cdot \begin{pmatrix} D_{yz}g \\ -D_{xz}g \end{pmatrix}  \neq 0.
\end{eqnarray*}

This provides the full unfolding of the cusp under the independent variation of two slow parameters. 

\begin{remark}
These conditions should be compared to the defining equations in \cite{broer2013}, and in particular the transversality condition, which is a corrected version of the nondegeneracy condition (A) in their paper. We note that the appropriate transversality condition is correctly identified later in equation (35) of \cite{broer2013}, after a series of coordinate transformations. We also refer the reader to a blow-up analysis of the cusp singularity for standard systems in \cite{jardon2016}. 
\end{remark}

\subsubsection{Versal slow unfoldings of $A_c$ singularities in the standard form}
\eIL{
Versal slow unfoldings of higher-order $A_c$ singularities in the standard case can be characterised explicitly by using classical formulas as shown in the proof of Lemma \ref{lemma:arbitrary}. Slow-fast systems which can locally be placed in the standard form

\begin{eqnarray*}
\begin{pmatrix} z_1'  \\ \vdots \\ z_{c}' \\ x' \end{pmatrix} &=& \begin{pmatrix} 0 \\ \vdots \\ 0 \\ 1 \end{pmatrix} (x^{c+1} + z_{c} x^{c-1} + \cdots + z_{2} x + z_{1}) + \eps G(z_1,\cdots,z_c,x,\eps)
\end{eqnarray*}

will exhibit an $A_c$ singularity of the critical manifold $\{(z_1, \cdots, z_{c-1},x) \in \mathbb{R}^c:  (x^{c+1} + z_{c} x^{c-1} + \cdots + z_{2} x + z_{1}) = 0\}$ with the fast fibers at the origin. Generally for fixed $c$, the corresponding $A_c$ singularity requires $c$ slow variables for a slow unfolding.}\\

%

\subsection{Three-component negative feedback oscillator} \label{sec:threecomp}
A fundamental characteristic of biochemical oscillators is the presence of negative feedback with time delay. Novak and Tyson considered several examples of biochemical networks, including autonomous systems and delay differential equations. They argued that sufficienly many intermediate steps can model the effect of a delay in a negative feedback loop, thus generating sustained oscillations   \cite{novak2008}. We consider the following minimal, autonomous three-component model studied in \cite{wechselberger2020}:

\begin{eqnarray}
\begin{pmatrix} x' \\ y' \\ z' \end{pmatrix} &=& \begin{pmatrix} \alpha_1 \left( \frac{1}{1+z^2} - x \right) \\ \alpha_2 x - 1 \\ \alpha_3 (y-z) \end{pmatrix} y + \eps \begin{pmatrix}  \alpha_1 \left( \frac{1}{1+z^2} - x \right) \\ \alpha_2 x \\ \alpha_3 (y-z) \end{pmatrix}. \label{eq:threecompfull}
\end{eqnarray}

The (dimensionless) parameters are $\alpha_1,\alpha_3 > 0$, $\alpha_2 > 1$, and $\eps \ll 1$. The system is in nonstandard form \eqref{eq:nonstandard}. The layer problem of \eqref{eq:threecompfull} is given by

\begin{eqnarray}
\begin{pmatrix} x' \\ y' \\ z' \end{pmatrix} &=& \begin{pmatrix} \alpha_1 \left( \frac{1}{1+z^2} - x \right) \\ \alpha_2 x - 1 \\ \alpha_3 (y-z) \end{pmatrix} y. \label{eq:threecomplayer}
\end{eqnarray}

Here we have $N(x,y,z)$ the vector function on the RHS and $f(x,y,z) = y$. The (regular part of the) critical manifold is given by the plane
\begin{eqnarray*}
S &=& \{y = 0\}.
\end{eqnarray*}

The nontrivial eigenvalue is given by the scalar function
\begin{eqnarray*}
D\!f   N|_S &=& \alpha_2 x - 1.
\end{eqnarray*}
The contact set is the curve
\begin{eqnarray*}
F &=& \{(x,y,z) \in S: x = 1/\alpha_2\}.
\end{eqnarray*}

\begin{figure}[!]
  \centering
  (a)      \includegraphics[height=0.5\textwidth,width=0.9\textwidth]{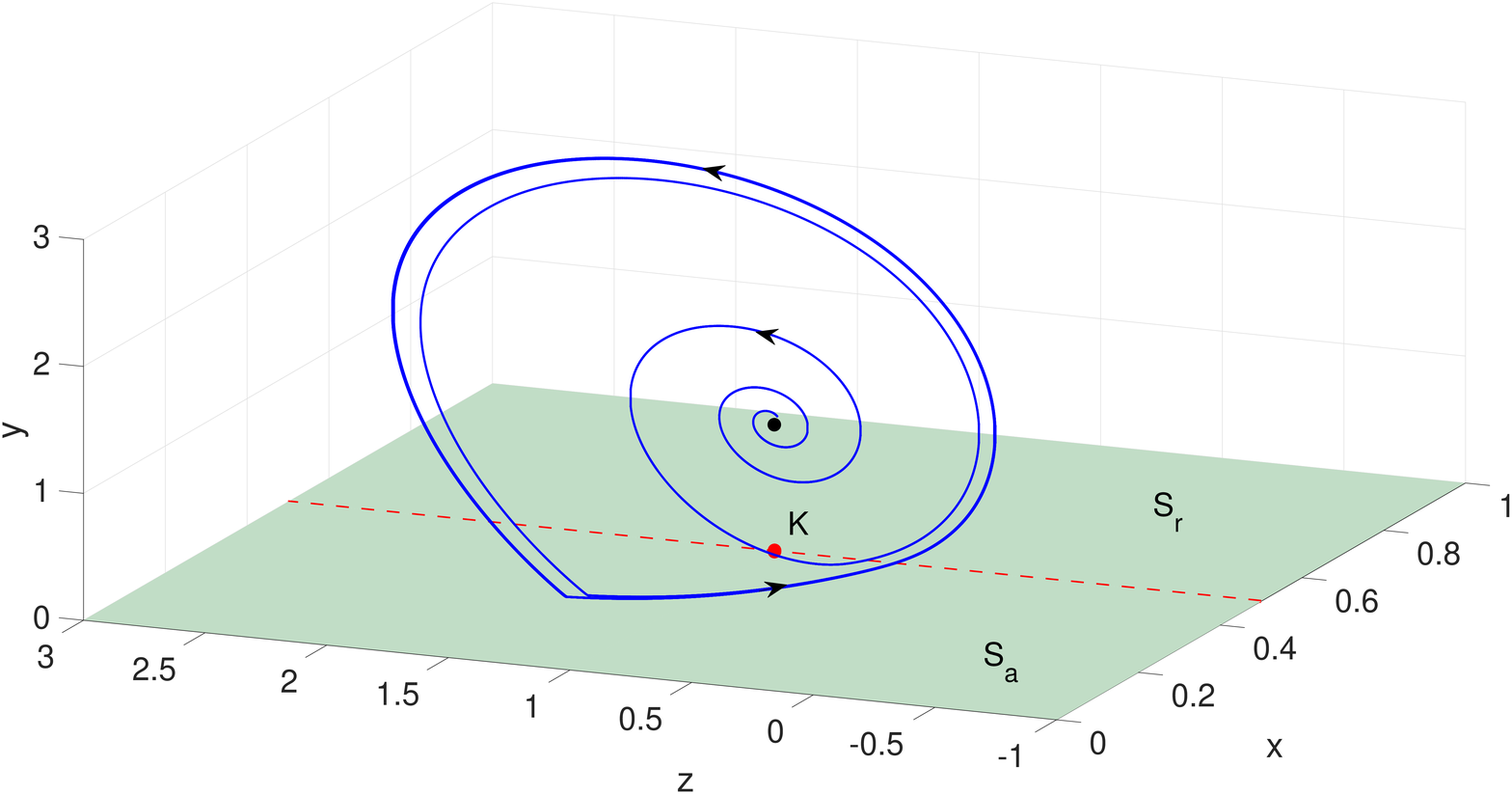}\\
     (b)           \includegraphics[height=0.5\textwidth,width=0.9\textwidth]{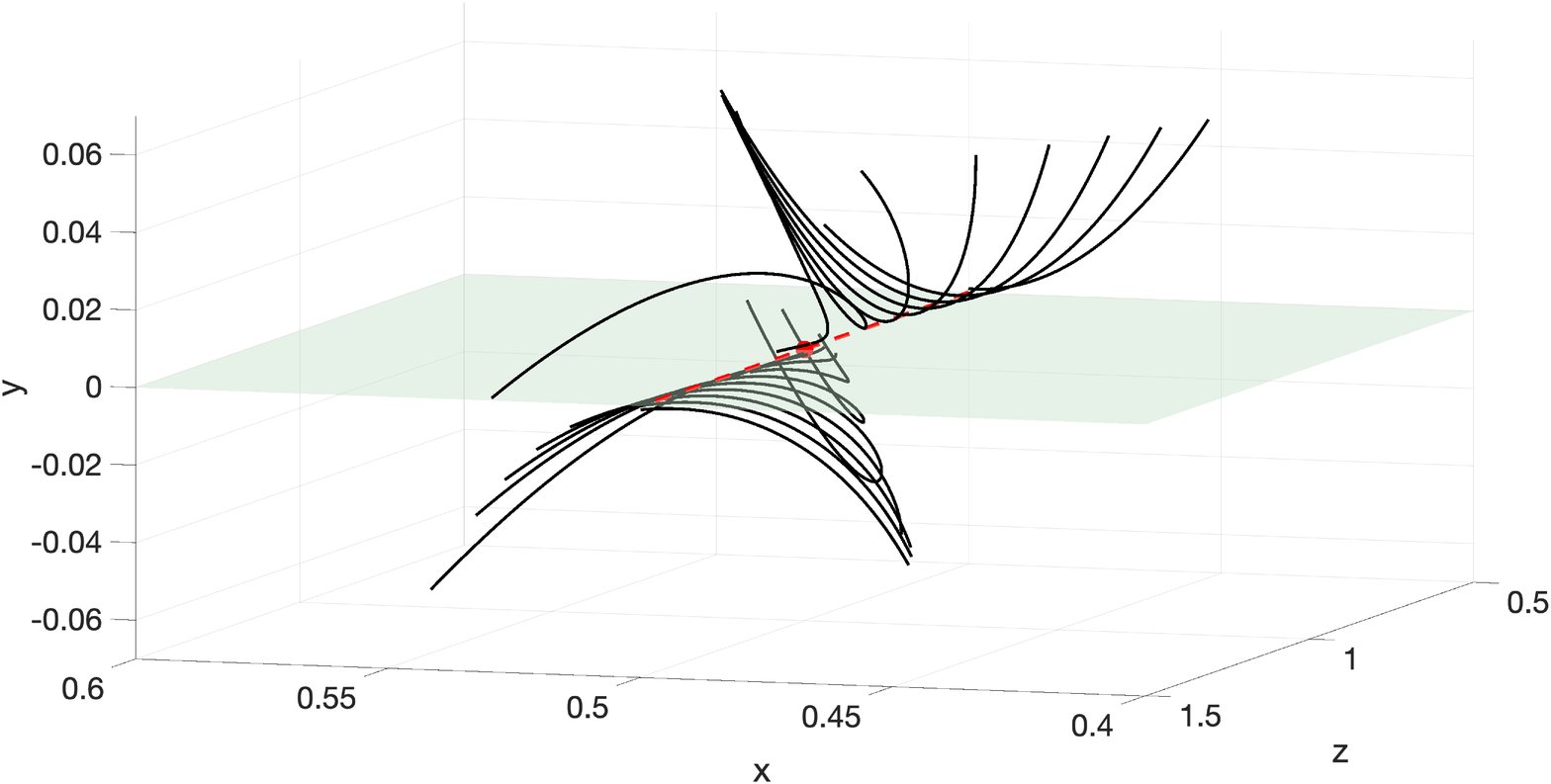}
      \caption{Sample trajectory (blue curve) of the three-component negative feedback oscillator \eqref{eq:threecompfull} with initial condition $(x,y,z) \approx (0.5198, 1.0205, 1.0205)$ and parameter set $(\eps,\alpha_1, \alpha_2, \alpha_3) = (0.0005,0.2,2,0.2)$. An isolated contact point $K$ (red point) of order 2 (\eqref{eq:negfeedbackcusp}) lies on the fold line (red dashed line). A saddle-focus equilibrium point (black point) of \eqref{eq:threecompfull} is also shown. (b) Local geometry of the fast fibers near the line of folds (note that this figure is rotated to better show the curvature).}  
      \label{fig:threecomp}
\end{figure}

The layer flow \eqref{eq:threecomplayer} also admits an isolated saddle-focus equilibrium point 
\begin{eqnarray*}
q &=& (1/\alpha_2, \sqrt{\alpha_2 - 1},\sqrt{\alpha_2 - 1}),
\end{eqnarray*}
given by the zero set $\{N(x,y,z) = 0\}$. This point persists as an equilibrium point for the full system  \eqref{eq:threecompfull} for small nonzero values of $\eps$. The saddle-focus is responsible for the bending of the fast fibers, generating a global return mechanism for the observed relaxation oscillations  (see Fig. \ref{fig:threecomp}). This periodic orbit can be decomposed into a slow segment near $S$ which crosses $F$, and a fast global reinjection arising from intersections between the two-dimensional unstable manifold of the saddle-focus and the fast fiber bundle near $S$. Further details on the global dynamics are provided in \cite{wechselberger2020}; for instance, the curve $F$ divides $S$ into an attracting and a repelling branch, denoted $S_a$ resp. $S_r$ in Fig. \ref{fig:threecomp}.   \\

{\bf Contact folds.} Almost all points on $F$ are fold points. We observe this by checking the conditions in Lemma \ref{lemma:fold}(a)-(b):

\begin{eqnarray*}
\text{rank}(D\!f  )|_{F} &=& n -k = 1\\
D\!f  D\!N N|_{F} &=& \alpha_1\frac{\alpha_2 - (1+z^2)}{1+z^2}.
\end{eqnarray*}

(observe how simple the nondegeneracy condition becomes in the codimension-one case: we have $\text{adj}(D\!f   N) = 1$). 

Thus, the parabolic coefficient is nontrivial everywhere on the contact set except where $z = \sqrt{\alpha_2 - 1}$. We call this distinguished point 
\begin{eqnarray}
K &=& \{(x,y,z) \in F: z = \sqrt{\alpha_2 - 1}\} = \{1/\alpha_2, 0, \sqrt{\alpha_2-1}\} .\label{eq:negfeedbackcusp}
\end{eqnarray}

At $K$, the conditions in Lemma \ref{lemma:cusp}(a)--(b) are satisfied. 

\begin{remark} \label{rem:kminus}
We note the existence of an unphysical contact point $K_-$ of order at least two, with the component $z = -\sqrt{\alpha_2 - 1}$. \\
\end{remark}

{\bf Contact cusps.} We verify that $K$ is a contact cusp. We check Lemma \ref{lemma:cusp}(c) The identities $D^2\!f = 0$ and $D^3\!f = 0$ greatly simplify the calculations; we need only evaluate $D\!f   D^2N (N,N,1)$ and $D\!f   D\!N D\!N N$. The second derivative in the first term admits the following simple formula in the codimension-one case:

\begin{eqnarray*}
D\!f   D^2(N,N,1) &=& D\!f  \begin{pmatrix} N^T H\!N_1 N \\  N^T H\!N_2 N  \\ \vdots \\ N^T H\!N_n N \end{pmatrix},
\end{eqnarray*}

where $H\!N_i$ denotes the Hessian of the scalar function $N_i (z)$. This term evaluates to 0, which can be read off from the fact that $N_2(z)$ is only linear, whereas $D_xf$ and $D_zf$ are both zero. On the other hand, the last term is nontrivial:

\begin{eqnarray*}
D\!f   D\!N D\!N N &=& \frac{2\alpha_1 \alpha_3 (\alpha_2 - 1)}{\alpha_2}.
\end{eqnarray*}

As long as $\alpha_2 > 1$, $K$ has contact-order of 2 (note that $K$ exists for $\alpha_2 \geq 1$).

We now check the remaining transversality condition, Lemma \ref{lemma:cusp}(d). As before, we write down

\begin{eqnarray*}
C_0(K) &=& \left.\begin{pmatrix} D\!f   \\ D^2\!f N + D\!f   D\!N  \end{pmatrix} \right\vert_{z_0}.
\end{eqnarray*}

at the contact  point $K$, we have

\begin{eqnarray*}
D\!f   &=& \begin{pmatrix} 0 & 1 & 0 \end{pmatrix}\\
D\!N &=& \begin{pmatrix} \alpha_1 & 0 & -2(\alpha_1/\alpha_2) \sqrt{\alpha_2-1}\\ \alpha_2 & 0 & 0 \\ 0 & \alpha_3 & -\alpha_3  \end{pmatrix},
\end{eqnarray*}

so

\begin{eqnarray*}
C_0(K) &=& \begin{pmatrix} 0 & 1 & 0 \\ \alpha_2 & 0 & 0  \end{pmatrix},
\end{eqnarray*}

which has the maximal rank of $2$ for the parameter values we consider. The three-component feedback oscillator therefore exhibits a contact cusp at the point $K$ (see \eqref{eq:negfeedbackcusp}).


\subsection{Mitotic oscillator} \label{sec:mitotic}

\begin{figure}[!]
  \centering
(a)        \includegraphics[height=0.5\textwidth,width=0.9\textwidth]{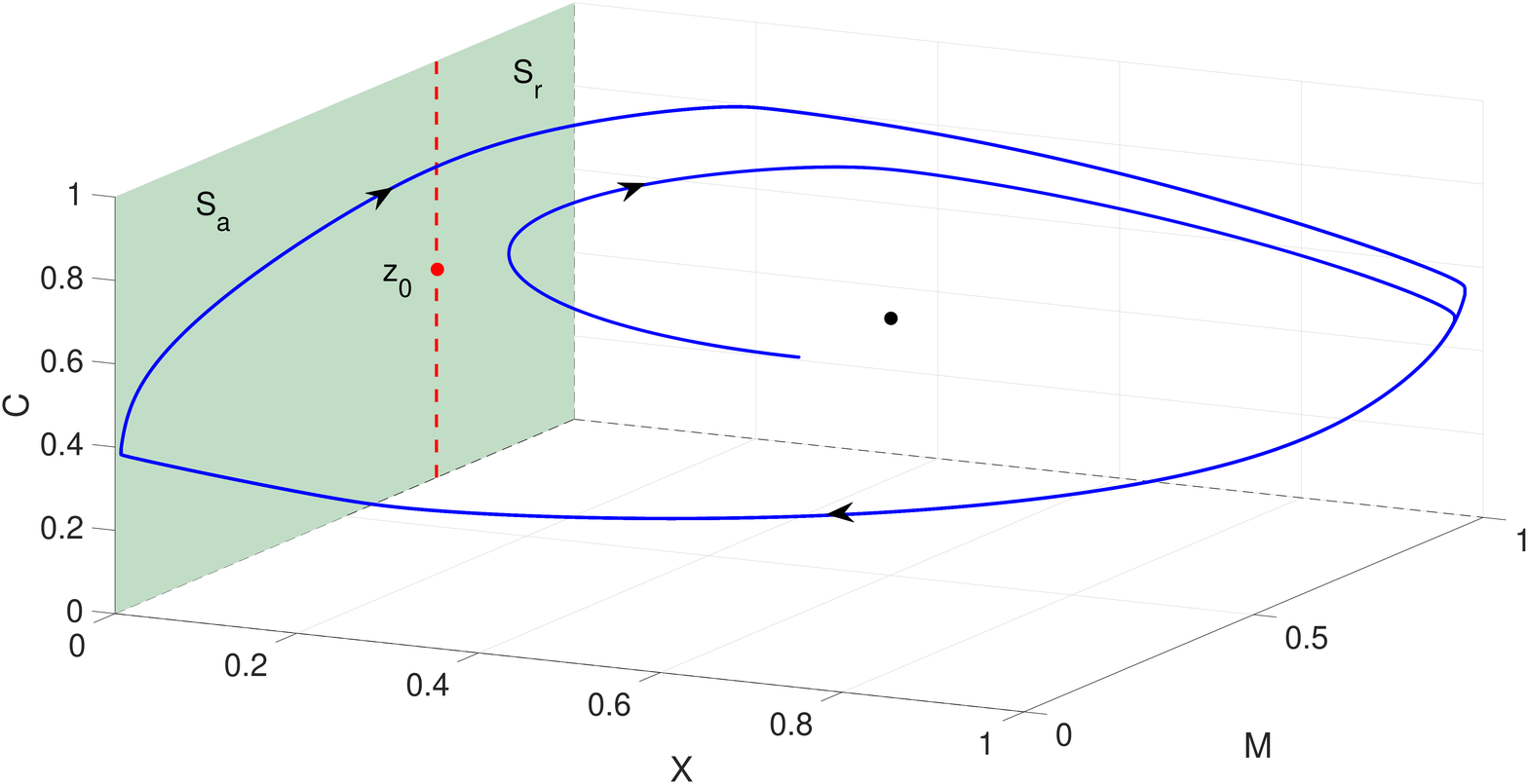}\\
(b)        \includegraphics[height=0.5\textwidth,width=0.9\textwidth]{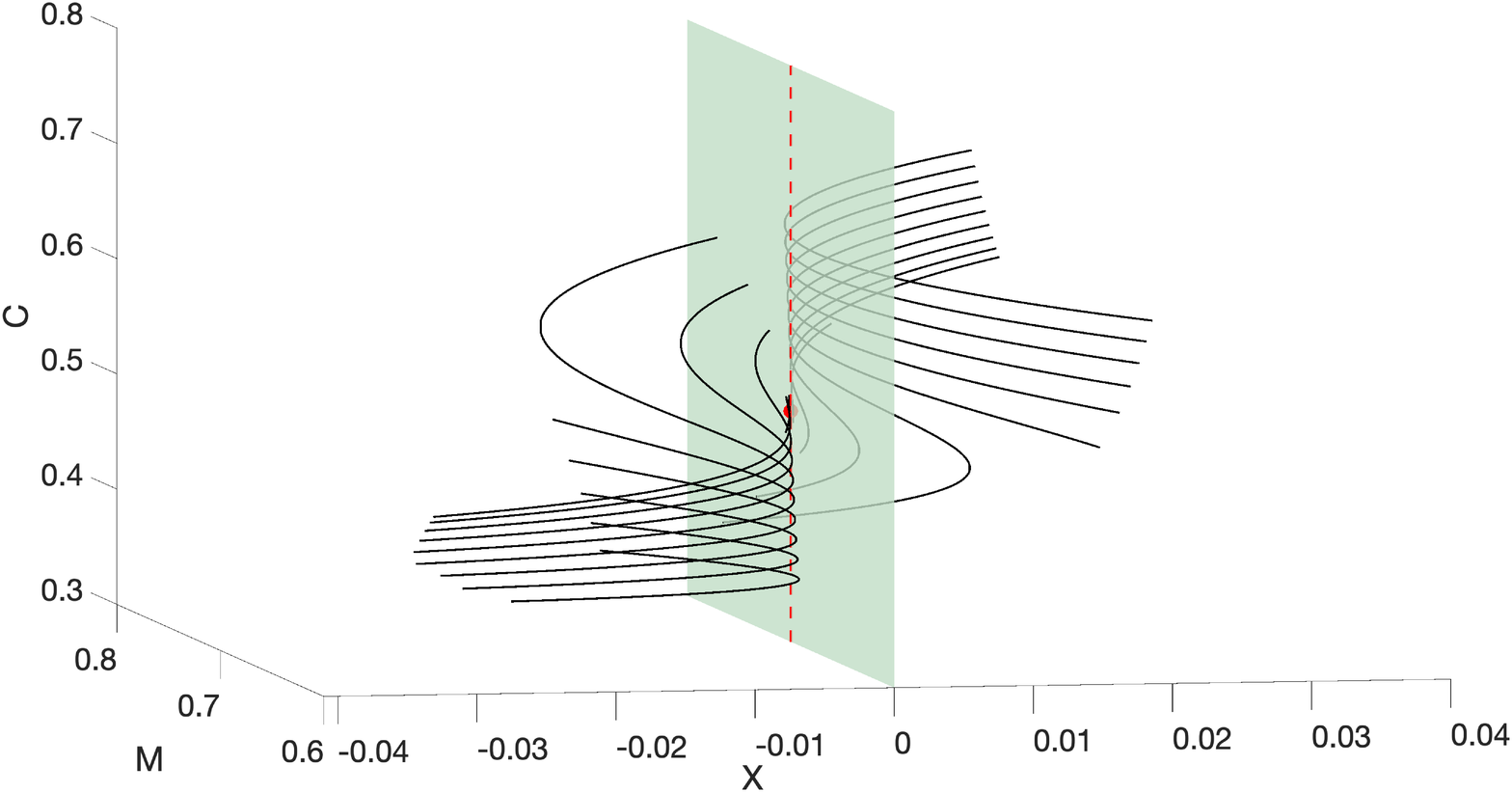}
      \caption{Sample trajectory (black curve) of the system \eqref{eq:mitotic} with initial condition $(X,M,C) = (0.5,0.5,0.5)$ and $\eps = 0.0021$. The green plane $X = 0$ lies in the set of equilibria of \eqref{eq:mitotic}.  An isolated contact point $K$ (red point) of order 2 (\eqref{eq:mitoticcusp}) lies on the fold line $F = \{X = 0\} \cap \{M = 7/10\}$ (red dashed line). An isolated saddle-focus equilibrium of \eqref{eq:mitotic} is also depicted (black point). (b) Local geometry of the fast fibers near the fold line. Note that we extend the fast fibers into the unphysical regime $X < 0$.}  
      \label{fig:mitotic}
\end{figure}

We demonstrate the existence of a cusp in Goldbeter's minimal model for the embryonic cell cycle \cite{goldbeter1991}. The original formulation contains terms of Michaelis-Menten type to study the existence of sustained oscillations due to negative feedback loops.  An analysis from the GSPT point of view is provided in \cite{kosiuk2015}, where an isolated, strongly attracting limit cycle is proven to exist for sufficiently small values of a singular perturbation parameter (see Fig. \ref{fig:mitotic}). Following their formulation, consider the system

\begin{eqnarray}
\frac{dX}{dt} &=& \left( M(1-X)(\eps + X) - \frac{7}{10} X(\eps + 1 - X) \right)  F_{\eps}(M) \nonumber\\
\frac{dM}{dt} &=& \left( \frac{6C}{1+2C}(1-M)(\eps + M) - \frac{3}{2}M(\eps + 1 - M) \right) F_{\eps}(X) \label{eq:mitotic}\\
\frac{dC}{dt} &=& \frac{1}{4}(1-X-C)F_{\eps}(X,M), \nonumber
\end{eqnarray}

where 

\begin{eqnarray*}
F_{\eps}(X,M) &=& F_{\eps}(X)F_{\eps}(M)\\
F_{\eps}(X) &=& (\eps + 1 - X) (\eps + X)\\
F_{\eps}(M) &=& (\eps + 1 - M)(\eps + M).
\end{eqnarray*}

The layer problem is given by
\begin{eqnarray*}
\frac{dX}{dt} &=& \left( M - \frac{7}{10} \right)F_0(X,M,C)\\
\frac{dM}{dt} &=& \left( \frac{6C}{1+2C}-\frac{3}{2} \right)F_0(X,M,C)\\
\frac{dC}{dt} &=& \frac{1}{4}(1-X-C)F_0(X,M,C),
\end{eqnarray*}

for 
\begin{eqnarray*}
F_0(X,M,C) &=& XM(1-X)(1-M).
\end{eqnarray*}

The critical manifold $S$ is given by regular two-dimensional subsets of the zero set $\{F_0(X,M,C) = 0\} = \{X = 0\} \cup \{X = 1\} \cup \{M = 0\} \cup \{M = 1\}$. The four faces intersect at four `corners,' and blow-up is necessary to analyze the dynamics on those lines (see \cite{kosiuk2015}).  \\
The one-dimensional linear fast fibers are spanned by the vector
\begin{eqnarray*}
N(X,M,C) &=& \begin{pmatrix}   M - \frac{7}{10} \\ \frac{6C}{1+2C}-\frac{3}{2}  \\ \frac{1}{4}(1-X-C) \end{pmatrix}
\end{eqnarray*}

at points $(X,M,C) \in S$. In the sequel we denote $f = F_0$ so that we can read off the defining equations classifying the singularities along the contact set. \\

Let us record the following derivatives:
\begin{eqnarray*}
D\!f   &=& \begin{pmatrix} M(1-M)(1-2X) & X(1-X)(1-2M) & 0 \end{pmatrix}\\
D^2\!f &=& \begin{pmatrix} 2M(M-1) & (1-2X)(1-2M) & 0 \\(1-2X)(1-2M) & 2X(X-1) & 0 \\ 0 & 0 & 0  \end{pmatrix}\\
D\!N &=& \begin{pmatrix} 0 & 1 & 0 \\ 0 & 0 & \frac{6}{(2C+1)^2} \\ -\frac{1}{4} & 0 & -\frac{1}{4} \end{pmatrix}
\end{eqnarray*}

We now restrict ourselves to the plane $S = \{X = 0\}$ (see Fig. \ref{fig:mitotic}). We have
\begin{eqnarray*}
D\!f   N|_{S} &=& M(1-M)(M-7/10).
\end{eqnarray*}

The critical manifold $S$ loses normal hyperbolicity along the lines $M = 0$, $M = 7/10$, and $M = 1$, and $S$ is attracting on the subset $S_a = S \cap \{0 < M < 7/10\}$ and repelling on the subset $S_r = S \cap \{7/10 < M < 1\}$. We do not consider the degenerate lines where $S$ intersects the faces $M = 0$ and $M = 1$. We focus on the fold line $F = S \cap \{M = 7/10\}$. Note that the critical manifold remains locally two-dimensional along this line, but the matrix $D\!f   N$ drops rank along $F$. The left- and right-nullvectors of $D\!f   N|_F$ are $l = r = 1$ .\\

{\bf Contact folds.} We test the nondegeneracy condition Lemma \ref{lemma:fold}(b).

\begin{eqnarray*}
l (D^2\!f(Nr,Nr) + D\!f  D\!N(Nr,r)) &=& N^T D^2\!f N + D\!f  D\!NN\\
&=& 0 + \frac{63}{200}\frac{2C-1}{2C+1}
\end{eqnarray*}

along $F$. Thus, a line of fold points separates $S_a$ from $S_r$, but there is a distinguished point 
\begin{eqnarray}
K &=& F \cap \{C = 1/2\} = \{(X,M,C) = (0,7/10,1/2)\}, \label{eq:mitoticcusp}
\end{eqnarray}

which has higher contact order. At $K$, the conditions Lemma \ref{lemma:cusp}(a)--(b) are satisfied.\\

{\bf Contact cusps.} We test the nondegeneracy condition Lemma \ref{lemma:cusp}(c):

\begin{eqnarray*}
&& l\cdot ( D^3\!f(Nr,Nr,Nr) + 3D^2\!f (D\!N(Nr,r),Nr)  + D\!f   D^2N(Nr,Nr,r) + D\!f   D\!N(D\!N(Nr,r),r))\\
&=&D^3\!f(N,N,N) + 3D^2\!f (D\!N N,N)  + D\!f   D^2N(N,N,1) + D\!f   D\!N D\!N N\\
&=&  0 + 0 + 0 + \frac{63}{1000}
\end{eqnarray*}

when evaluated at $K$.\\

We test the transversality condition Lemma \ref{lemma:cusp}(d):

\begin{eqnarray*}
C_0(K) &=& \left.\begin{pmatrix} D\!f   \\ N^T D^2\!f  + D\!f   D\!N  \end{pmatrix} \right\vert_{K}.
\end{eqnarray*}

At $K$ we have

\begin{eqnarray*}
D\!f   &=& \begin{pmatrix} 21/100 & 0 & 0 \end{pmatrix}\\
D^2\!f &=& \begin{pmatrix} -21/50 & -2/5 & 0 \\ -2/5 & 0 & 0 \\ 0 & 0 & 0 \end{pmatrix}\\
N &=& \begin{pmatrix} 0 \\ 0 \\ 1/8 \end{pmatrix}\\
D\!N &=& \begin{pmatrix} 0 & 1 & 0 \\ 0 & 0 & 3/2 \\ -1/4 & 0 & -1/4 \end{pmatrix},
\end{eqnarray*}

and so altogether we have

\begin{eqnarray*}
C_0(K) &=& \begin{pmatrix} 21/100 & 0 & 0 \\ 0 & 21/100 & 0 \end{pmatrix}.
\end{eqnarray*}

The mitotic oscillator therefore exhibits a contact cusp at the point $K$ (see \eqref{eq:mitoticcusp}). 

\begin{remark}
This analysis may repeated for the other three faces $\{X = 1\}$, $\{M=0\}$, and $\{M = 1\}$. There are three additional lines of contact order at least one:
\begin{eqnarray*}
L_{X=1} &=& \{X=1\} \cap \{M = 7/10\}\\
L_{M=0} &=& \{M=0\} \cap \{C = 1/2\}\\
L_{M=1} &=& \{M=1\} \cap \{C = 1/2\}.
\end{eqnarray*}

Away from the corners where the faces intersect (i.e. for the ranges $0 < X < 1$ and $0 < M < 1$), the points on these lines are all contact folds except for the following three contact cusps:
\begin{eqnarray*}
K_{X=1} &=& (1,7/10,1/2)\\
K_{M=0} &=& (1/2,0,1/2)\\
K_{M=1} &=& (1/2,1,1/2).
\end{eqnarray*}
\end{remark}

\section{Concluding remarks} \label{sec:conclusion}

We have given a rigorous classification of the contact singularities of singularly perturbed systems in the nonstandard form \eqref{eq:nonstandard}, and \eIL{we provided computable coordinate-independent criteria to identify slow-generic contact singularities of arbitary order (Lemmas \ref{lemma:fold} -- \ref{lemma:arbitrary})}. We emphasize that the correct coordinate-free conditions for the corresponding slow unfoldings are {\it not} immediately given by the usual defining equations from singularity theory: in the classical context, the parameter directions are already assumed to have been located, whereas the slow variables are generally not explicitly identified in the system \eqref{eq:nonstandard}. \\

Constructing a vector field factorisation \eqref{factor} for systems in the general form \eqref{eq:nonstandard} is a nontrivial first step. We expect that the procedure of first applying constructive factorisation algorithms \cite{goeke2014}  and then identifying and classifying the points in the contact set can be automated in a wide variety of applications. We also expect that the defining equations can serve as test functions for continuation and bifurcation software like AUTO. \\

Finally, we point out that loss of normal hyperbolicity may occur by other means which are dynamically relevant. The contact set (Def. \ref{def:contactset}) can be defined more generally to include rank drops larger than one, corresponding to more degenerate contact scenarios between the critical manifold and surfaces in the fast fiber bundle. 
Away from the contact set, there may also be subsets of the critical manifold where complex-conjugate pairs of eigenvalues lie on the imaginary axis, which is associated with the existence of delayed Hopf \cite{baer1989,neishtadt1,neishtadt2} and singular Hopf bifurcations \cite{baer1986a,baer1986b}.  This paper thus serves as a starting point for a more complete theory of the loss of normal hyperbolicity for nonstandard slow-fast systems.

%
%
%
%

\newpage
\thispagestyle{plain}
\appendix \section{Appendix} \label{appendix}

We write down basic definitions for jet spaces and the contact group of diffeomorphisms (see for eg. \cite{golubitsky1973,izumiya2015} for a full treatment of the standard singularity theory).
\begin{definition} The $k$-jet space $J^{k}(n,m)$ of smooth germs $f:\mathbb{R}^n \to \mathbb{R}^m$ is defined by
\begin{eqnarray*}
J^{k}(n,m) &=& \mathscr{M}_n \cdot \mathscr{E}(n,m)/ \mathscr{M}^{k+1}_n \cdot \mathscr{E}(n,m),
\end{eqnarray*}

where
\begin{eqnarray*}
\mathscr{E}(n,m) &=& (\mathscr{E}_n)^m
\end{eqnarray*}

is the direct product of $m$ copies of the set $\mathscr{E}_n$ of smooth germs from $\mathbb{R}^n$ to $\mathbb{R}$,
\begin{eqnarray*}
\mathscr{M}_n &=& \mathscr{E}_n \cdot \{x_1, \cdots, x_n\}
\end{eqnarray*}

is the unique maximal ideal of germs vanishing at the origin, and
\begin{eqnarray*}
\mathscr{M}^k_n &=&  \mathscr{E}_n \cdot \{ x_1^{i_1}, \cdots, x_n^{i_n}, i_1 + \cdots + i_n  = k\}
\end{eqnarray*}

is the set of germs with vanishing partial derivatives of order less than or equal to $k-1$ at the origin.
\end{definition}

\begin{remark}The set $J^k(n,m)$ may be identified with the set of polynomials of total degree less than or equal to $k$. \end{remark}

The definition of contact classes used in the paper is due to Mather:\\

\begin{definition} \label{def:contactequiv} The {\it contact group} $\mathscr{K}$ is the set of germs of diffeomorphisms of $(\mathbb{R}^n \times \mathbb{R}^m, (0,0))$ which can be written in the form

\begin{eqnarray*}
H(x,y) &=& (h(x),H_1(x,y)),
\end{eqnarray*}

where $h$ acts on the right (i.e. $h \cdot f = f \circ h^{-1}$) and $H_1(x,0) = 0$ for $x$ near 0.  We say that $f$ is {\it $\mathscr{K}$-equivalent} to $g$ if $g$ lies in the group orbit of $f$.  We refer to this as the contact class of $f$. \end{definition}

\begin{remark} Suppose $f,g \in \mathscr{M}_n \cdot \mathscr{E}(n,m)$ and $k = (h,H) \in \mathscr{K}$. Then $g = k \cdot f$ if and only if
\begin{eqnarray*}
(x,g(x)) &=& H(h^{-1}(x),f(h^{-1}(x))).	
\end{eqnarray*}
Observe that $H$ sends the graph of $f$ to the graph of $g$ near 0 (i.e. the zero sets of $\mathscr{K}$-equivalent germs are diffeomorphic). \end{remark}

\end{document}